\documentclass[12pt]{amsart}

\newtheorem{theorem}{Theorem}[section]
\newtheorem{lemma}[theorem]{Lemma}
\newtheorem{proposition}[theorem]{Proposition}

\newtheorem{example}[theorem]{Example}
\newtheorem{definition}[theorem]{Definition}
\newtheorem{corollary}[theorem]{Corollary}

\renewcommand{\theequation}
{\arabic{equation}}
\newcommand{\Co}{\mbox{$\mathbb{C}$}}
\newcommand{\N}{\mbox{$\mathbb{N}$}}
\newcommand{\T}{\mbox{$\mathbb{T}$}}
\newcommand{\Z}{\mbox{$\mathbb{Z}$}}
\newcommand{\K}{\mbox{$\mathbb{K}$}}
\newcommand{\Di}{\mbox{$\mathbb{D}$}}
\newcommand{\Rl}{\mbox{$\mathbb{R}$}}

\newcommand{\C}{\mbox{${\mathcal C}$}}
\newcommand{\D}{\mbox{${\mathcal D}$}}

\newcommand{\G}{\mbox{${\mathcal G}$}}
\newcommand{\R}{\mbox{${\mathcal R}$}}

\newcommand{\M}{\mbox{${\mathcal M}$}}

\begin{document}
\title[Module equivalence for function algebras]{}
\title[Function modules]{Isomorphisms of function modules, \\
and generalized approximation in modulus}
\author{David Blecher}
\address{Department of Mathematics\\
University of Houston\\
Houston, TX 77204-3476 }
\email{dblecher@math.uh.edu}
\author{Krzysztof Jarosz}
\address{Department of Mathematics and Statistics \\
Southern Illinois University \\
Edwardsville, IL 62026-1653}
\email{kjarosz@siue.edu}
\date{Sept 4, 1999}

\begin{abstract}
For a function algebra $A$ we investigate relations between the following
three topics: isomorphisms of singly generated $A$-modules, Morita
equivalence bimodules, and `real harmonic functions' with respect to $A$.
We also consider certain groups which are naturally associated
with a uniform algebra $A$.
We illustrate the notions considered with several examples.
\end{abstract}

\maketitle

\section{Introduction.}

By a {\em uniform algebra} or {\em function algebra} 
on a compact Hausdorff space $\Omega$, we shall mean 
a subalgebra $A$ of $C(\Omega)$ (the continuous complex valued functions
on $\Omega$) which
contains constants
and separates points.  In most
of this paper we are concerned with closed submodules of $C(\Omega)$ of the 
form $Af$, where $f$ is a strictly positive and continuous function on 
$\Omega$.  In Part C we will allow $f$ to be nonnegative.

Before we proceed any further with this introduction, we shall 
take a paragraph to explain why these modules are more general 
than they appear to be at first.  
We need a little notation.  By a 
{\em concrete function $A$-module} we shall mean a 
closed linear subspace
of $C(K)$, for a compact Hausdorff space $K$, which is
closed under multiplication by $\pi(A)$, where
$\pi : A \rightarrow
C(K)$ is a unital homomorphism.  By an {\em 
(abstract) function $A$-module} we shall mean a
Banach $A$-module $X$
which is isometrically $A$-module isomorphic to  
(in future we shall simply say `$A$-isometric to' for short) a
concrete function module.   Although we shall not particularly use 
this here, this class of modules was given several 
equivalent abstract characterizations in \cite{BLM}; for example it 
coincides with the class of Banach $A$-modules whose module action
is contractive with respect to the injective tensor product.
It is proved in \cite{Bsh}, that any algebraically singly generated 
 faithful\footnote{A left module
$X$ is faithful if $aX = 0$ implies $a = 0$.} function $A$-module
is $A$-isometric to one of the form $Af$ described in the 
first paragraph, for some strictly positive 
continuous function $f$ on some (possibly different) 
compact $\Omega$.

In this paper we investigate the relations
between the following three topics: isometries and almost isometries 
between modules of the type discussed above, real `harmonic' functions 
on $\Omega $ with respect to $A$,
and Morita equivalence bimodules over $A$.

In Part A we provide a necessary and sufficient condition on
functions $f_{1},f_{2}$ for the
modules $Af_{1},Af_{2}$ to be isometric (as
Banach spaces). The result generalizes the
 classical description of isometries
of uniform algebras. We then extend the result to almost isometric 
modules; where modules $Af_{1},Af_{2}$ are called `almost isometric' if 
for any $\varepsilon >0$ there is a surjective linear isomorphism $%
T:Af_{1}\rightarrow Af_{2}$ such that $\left\| T\right\| \left\|
T^{-1}\right\| <1+\varepsilon $.

One of the main results in Part A, is a characterization of 
modules of the type $Af$ which are almost isometric to $A$.  
If $\Omega = \partial A$, the Shilov boundary of $A$, we will see that 
``$Af \cong A$ almost isometrically'' is equivalent to 
$f$ being uniformly
approximable by the moduli of invertible elements of $A$.
That is, $f$ is in the
 uniform closure $\bar{Q}$ of $Q$, where $Q = \{ |a| : a \in A^{-1} \}$.
Indeed, up to $A$-isometric isomorphism,
for any $\Omega$ on which $A$ is a function algebra,
the Banach $A$-modules which are almost
$A$-isometric to $A$, are exactly 
the submodules
$Af \subset C(\Omega)$, for some strictly positive
$f \in \bar{Q}$.  
The set $H_A(\Omega) = \{ \log f : 0 < f \in \bar{Q} \}$ is 
well known to those familiar with the theory of 
uniform algebras.  In 
particular, $A$ is known as `logmodular' if $H_A(\Omega) =
C_{\Rl}(\Omega)$.  For any function algebra $A$, the 
class $H_A(\Omega)$ 
deserves to be called a `harmonic class' of functions with 
respect to $A$.
In this paper, to be more specific, by a `harmonic class'  we 
shall mean a class $\mathcal{B}(\Omega )$
of real continuous functions on $\Omega ,$ which 
have at least the following properties:

\begin{itemize}
\item[(i)]  If $f\in \mathcal{B}(\Omega )$ then $f_{|_{\partial A}}\in
\mathcal{B}(\partial A)$.

\item[(ii)]  If $f\in \mathcal{B}(\Omega )$ then there exists a unique $%
\tilde{f}\in \mathcal{B}(M_{A})$ such that $\tilde{f}_{|_{\Omega }}=f$.  
($M_A$ is the maximal ideal space of $A$).

\item[(iii)]  Every $f\in \mathcal{B}(\Omega )$ achieves its maximum and
minimum value on $\partial A$.

\item[(iv)]  $f_{1},f_{2}\in \mathcal{B}(\Omega )$ implies $f_{1}+f_{2}\in
\mathcal{B}(\Omega )$ (Indeed most of the 
classes we study in this paper are additive groups).

\item[(v)]  If $A$ consists of functions which are analytic on a region $R$
in $\mathbb{C}^{n}$, and which separate points of $R$, then (via the obvious
homeomorphic embedding $R\rightarrow M_{A}$) the functions in $\mathcal{B}%
(M_{A})$ are genuinely harmonic on $R$.
\end{itemize}

This leads us conveniently into a description of Part B of the paper.
Here is a natural idea to attempt to generalize the class
$H_A(\Omega)$, or equivalently, the strictly positive functions 
which are uniform limits
 \begin{equation}
 f = \lim_n |k_n| \; ;  \; \; \text{where} \; k_n \; , \; 
h_n \;   \in A \;
\;  \text{with} \; 
 \; k_n h_n = 1  ,  \label{ac1}
\end{equation} 
on $\Omega$.   
We will write $A^{(n)}$ for the space of $n$-tuples with entries in $A$. An
element of $A^{(n)}$ will be called an $A$-tuple, and will often be regarded
as a function $\Omega \rightarrow \mathbb{C}^{n}$.  For two $A$-tuples
$H = (h_i), K = (k_i)$ of the same `length' we 
define $H.K = \sum_i h_i k_i \in A$.
Consider the set of strictly positive functions which are uniform limits
 \begin{equation}
f(w) = \lim_n \Vert K_n(w) \Vert_2 \; ;
  \; \; \text{where} \; K_n \; , \; H_n \; \in A^{(m_n)}
 \; \;  \text{with} \; \; K_n . H_n = 1  , 
 \label{ac2}
\end{equation} 
on $\Omega$, where the
$m_n \in \N$.  This looks like a natural generalization of 
$(\ref{ac1})$,
however one quickly sees that there is a hidden condition in $(\ref{ac1})$
which is not in $(\ref{ac2})$, which results in $(\ref{ac2})$ not 
corresponding to a `harmonic class'.  Namely, in $(\ref{ac1})$, because
$h_n = k_n^{-1}$, we automatically have $|h_n| \rightarrow f^{-1}$ 
uniformly.  We therefore
define $\M_A(\Omega)$ to be the set of strictly positive
functions $f$ on $\Omega$ which satisfy $(\ref{ac2})$ and also:
\begin{equation}
f(w)^{-1} = \lim_n \Vert H_n(w) \Vert_2 \; \;    \label{ac3}
\end{equation}
uniformly on $\Omega$.   This  might 
loosely be called a `tight convex approximation in 
modulus'.    One finds that now $\log \M_A(\Omega)$ 
is a `harmonic class', which contains
$H_A(\Omega)$.   This is shown in Part C.

We show
that $f \in \M_A(\Omega)$ if and only if $Af$ is a strong Morita
equivalence $A-A$-bimodule, with `inverse bimodule' $Af^{-1}$.
The notion of
strong Morita equivalence was defined and studied in \cite{BMP}, but 
for function algebras and singly generated $A$-modules $X$, this
notion may be 
viewed as a generalization of the notion  $X \cong A$.  Indeed 
$X$ is a `rank one' strong Morita
equivalence bimodule if and only if $X \cong A$ almost 
$A$-isometrically.    
We display certain groups that are naturally associated
with a uniform algebra $A$, such as the Picard group.
We also illustrate the notions considered with 
some interesting examples.

In Part C, we generalize still further.  We now allow
topologically singly generated function modules.   These 
correspond to submodules $(Af)^{\bar{}}$, the closure
taken in $C(\Omega)$, where $f$
is now allowed to be nonnegative continuous function on $\Omega$.
We define a larger
 harmonic class $\log \R(\Omega)$, 
which contains the harmonic classes mentioned 
earlier.  Just as the $\M$ class
corresponds to Morita equivalence, the $\R$ class corresponds
to the more general notion of `rigged module'.  Rigged modules
were intended to be a generalization of the notion of
`Hilbert C$^*-$module' and were studied in \cite{Bhmo,BMP};
but in our (singly generated) situation these are 
the modules $X$ for which the identity map $X \rightarrow X$
factors asymptotically, via contractive $A$-module maps, through
the `free' $A$-modules $A^{(n)}$.   
This asymptotic factorization may
be viewed as another generalization of the statement $X \cong A$ 
$A$-isometrically.   

 We show that for $A = A(\Di)$,
the disk algebra,
the topologically singly generated
 $A$-rigged modules are exactly (up to $A$-isometric isomorphism)
the modules of the form
$(Af)^{\bar{}}$, where $f$ is a continuous function
on $\bar{\Di}$ such that
 $f = |\phi|$ for some outer function $\phi \in H^\infty$.
This corresponds to the continuous nonnegative
functions on $\T$ whose logarithm is integrable.

We admit that one of the purposes of Parts 
B and C, was to begin to
 illuminate the 
function algebra case
of  the theory of  
 Morita
equivalence and rigged modules,
which the first author
 and various coauthors have
developed over the years   (see
\cite{Bnat} for a leisurely survey),
 and in particular 
to see the connections with  some problems 
concerning  function algebras.   Conversely, what we do here 
may lead to progress in the noncommutative situation.     
 
\section{Notations and definitions}

For a compact set $\Omega $ we denote by $C\left( \Omega \right) $ 
(resp. $C_{%
\mathbb{R}}\left( \Omega \right) $) the space of all complex valued (resp.
real valued) continuous functions on $\Omega $. 
For a function algebra $A$ on $\Omega$
we will write $M_{A}$ for the maximal ideal
space of $A$, and $\partial A$ for the Shilov boundary of $A$. Then $%
C(\partial A),C(M_{A}),$ and $C(\Omega )$, may be regarded, respectively, as
the minimal, maximal commutative, and generic commutative, C$^{\ast }-$%
algebra generated by $A$. $A$ may be viewed as a closed subalgebra of
continuous functions on either of these three compact spaces. For a set of
functions $\mathcal{E}$ we will write $\mathcal{E}_{+}$ for the nonnegative
functions in $\mathcal{E}$, and $\mathcal{E}^{+}$ for the strictly positive
functions in $\mathcal{E}$.  We will refer very often to the following 
important
subsets of $C(\Omega )$ which may be associated with $A$:

\begin{itemize}
\item  $A^{-1}$ - the set of invertible elements of $A$,

\item  $P=\{f\in C(\Omega ):f=|g|,g\in A\}$,

\item  $Q=\{f\in C(\Omega ):f=|g|,$ where $g\in A^{-1}\}$,

\item  $\mathcal{F}=\{f=\sum_{i=1}^{n}|g_{i}|^{2}\in C(\Omega ):n\in \mathbb{%
N},g_{1},\cdots g_{n}\in A\}$,

\item  $\mathcal{G}=\bar{\mathcal{F}}$, the uniform closure of 
$\mathcal{F}$ in $%
C(\Omega )$,

\item  $G_{A}(\Omega )$ is the closure of $ReA$ in 
$C_{\Rl}(\Omega )$,

\item  $H_{A}(\Omega )$ is the closure of $\{\log |a|:a\in A^{-1}\}$ in $%
C_{\Rl}(\Omega )$.
\end{itemize}

If we wish to specify the dependence on $\Omega$, we will write, 
for example, $Q(\Omega)$.  
It is well known that $G_{A}(\Omega ) \subset H_{A}(\Omega )$, and they 
are both harmonic classes in the sense of the introduction (a proof is 
also contained in Part C).
We recall that $A$ is called a `Dirichlet algebra' 
(resp. `logmodular algebra')
if $G_{A}(\partial A)=C_{\mathbb{R}}(\partial A)$ (resp. $H_{A}(\partial
A)=C_{\mathbb{R}}(\partial A)$).   The disk algebra $A(\Di)$ is
Dirichlet (and consequently logmodular).   Indeed this
is exactly saying
that the ordinary  Dirichlet problem (of harmonic extension
from the boundary) can be solved on the circle
$\T = \partial A$. 
This
of course, was Gleason's original reason for the name
`Dirichlet algebra' (see \cite{Gl}).   

We write $\hat{}:A\rightarrow C(M_{A})$ for
the Gelfand transform, and we often use the same symbol to denote an element
of an algebra and the corresponding Gelfand transform. There is a surprising
equivalence relation on $M_{A}$: $\phi ,\psi \in M_{A}$ are equivalent if
and only if $\Vert \phi -\psi \Vert <2$. The distinct equivalence classes
are called the (Gleason) parts \cite{Gl}. They are $\sigma -$compact subsets
of $M_{A}$. We refer the reader to \cite{sto} or \cite{Gam} as general
references on function algebras.

If $X_{1}$ and $X_{2}$ are Banach $A$-modules then we write $X_{1}\cong
X_{2} $ $A$-isometrically (resp. almost $A$-isometrically), if they are
isometrically $A$-isomorphic (resp. almost isometrically $A$-isomorphic).
As we said earlier, the latter term means that for any $\epsilon > 0$,
there exists an $A$-module isomorphism $T : X_1 \rightarrow X_2$
with $\Vert T \Vert \Vert T^{-1} \Vert < 1 + \epsilon$.  If, in the
above,  we replace 
the requirement that $T$ be an $A$-module map, with it being 
linear, then we simply say 
that $X_1 \cong X_2$ almost isometrically.

Whenever we use the words `singly
generated', it will be assumed, unless otherwise
qualified, to have the topological connotation.
Thus a `singly generated' Banach module $X$ over
$A$ has an element $x$ such that the closure
 of $Ax$ is $X$.   An {\em algebraically} singly generated
module $X$ has $Ax = X$ for some $x \in X$.
  
We will write $A^{(n)}$ for the space of $n$-tuples with entries in $A$. An
element of $A^{(n)}$ will be called an $A$-tuple, and will often be regarded
as a function $\Omega \rightarrow \mathbb{C}^{n}$.

For a Banach space $X$ we let $X_1$ be the unit ball of $X$, 
and we denote by $ext(X_1)$ the set of extreme points of the
unit ball of $X$. If $A$ is a subspace of $C\left( \Omega \right) $, and $x$
is a point in $\Omega $ we use the same symbol $x$ to denote the
corresponding functional on $A$ - namely, evaluation at $x$. There are several
definitions of the Choquet boundary $Ch(A)$ of a linear subspace
$A$ of $C\left( \Omega \right) $, here we adopt the following one:
\begin{equation}
Ch(A)=\left\{ x\in \Omega :x\in ext(A^{\ast}_1)\right\} .  \label{ChA}
\end{equation}
For general information, 
elementary Choquet theory  shows that this is the same as the set
$$\{ x \in \Omega : \delta_x \; \text{is the
only probability measure on} \; \Omega \; \text{extending} \;
(\delta_x)_{|_A} \; \}
$$
where $\delta_x$ is the Dirac mass considered as a functional on 
$C(\Omega)$.   We shall not use this explicitly here.
The Shilov boundary $\partial A$ is the closure of $Ch(A)$. As it may be
necessary to distinguish between functions from $A$ and its restrictions to $%
\partial A$, we will denote by $A_{|\partial A}$ the set of these
restrictions.

In the latter part of the paper,
we will be working with operator spaces.
However, usually, issues of `complete boundedness' do not
arise.  This is because for
 a linear operator $T$ mapping into
a subspace of a commutative
C$^*-$algebras we have $\Vert T \Vert = \Vert T \Vert_{cb}$.

\vspace{3 mm}

\begin{center}
Part A.
\end{center}

\section{Module isomorphisms}

It is well known \cite{N} that two uniform algebras are linearly isometric,
that is isometric as Banach spaces, if and only they are isomorphic as
algebras. In this section we show that similar results hold for function
modules of the form $Af$. To explain the main idea, let $A$ be a uniform
algebra on $\Omega$,
 and let $f_{1},f_{2}$ be strictly positive continuous
functions on $\Omega$. Suppose that 
there is an invertible function $g\in A$ such
that $\left| g\right| =\frac{f_{1}}{f_{2}},$ then
\begin{equation*}
Af_{1}\ni af_{1}\longmapsto agf_{2}\in Af_{2}
\end{equation*}
is an $A$-module isometry between $Af_{1}$ and $Af_{2}$. Suppose now that we
also have a homeomorphism $\varphi$ of $\Omega$ onto itself such that $%
\left\{ a\circ\varphi:a\in A\right\} =A$, and $\frac{f_{1}\circ\varphi }{%
f_{2}}=\left| h\right| $, for some invertible element $h$ of $A$. Then
\begin{equation}
Af_{1}\ni af_{1}\longmapsto (a\circ\varphi) hf_{2}\in Af_{2},  \label{can}
\end{equation}
is a linear isometry between $Af_{1}$ and $Af_{2}$, but this time it is not
an $A$-module isometry unless $\varphi$ is the identity map. We shall show
that any linear isometry between modules $Af_{1}$ and $Af_{2}$ is essentially
of the form (\ref{can}).  

We first need to introduce some  technical results.

\begin{lemma}
\label{ext}Suppose that $A$ is a subspace of $C\left( \Omega\right) $
separating the points of $\Omega$. Then

\begin{enumerate}
\item[\textbf{(i)}]  the original topology $\tau$ of $\Omega$ is identical
with the weak$^\ast$ topology $\sigma^{\ast}$ on $\Omega$ considered as a
subset of the unit ball of the dual space $A^{\ast}$,

\item[\textbf{(ii)}]  if $F$ is an extreme point of the unit ball of the
dual space $A^{\ast}$ then there is a point $x$ in $Ch\left( A\right) $ and
a scalar $\alpha$ of absolute value one such that $F=\alpha x$,

\item[\textbf{(iii)}]  for any $a\in A$ we have
\begin{equation}
\left\| a\right\| =\sup\left\{ \left| a\left( x\right) \right| :x\in
Ch(A)\right\} .  \label{KJ1}
\end{equation}
\end{enumerate}
\end{lemma}

\begin{proof}
\textbf{(i)} It is clear that the identity map $I:\left( \Omega,\tau\right)
\rightarrow\left( \Omega,\sigma^{\ast}\right) $ is continuous.
Since $\left(
\Omega,\tau\right) $ is compact it follows that the topologies are identical.

\textbf{(ii)} Assume that 
$F$ is an extreme point of the unit ball of the dual
space $A^{\ast}$. By the Krein-Milman Theorem,
 the set of norm one extensions
of $F$ to $C\left( \Omega\right) $ has an extreme point $\mu$.
It is easy to check that $\mu$ is also an extreme point of 
$C\left( \Omega\right)^{\ast}_1$ .
The argument is concluded by appealing to the  
very well-known fact that any extreme point of the last space
is of the form $\alpha x$ with $x\in \Omega,\left|
\alpha\right| =1.$

\textbf{(iii) }Let $a$ be a norm one element of $A$. By the Krein-Milman
Theorem there is an extreme point $F$ of the unit ball of $A^{\ast}$ such
that $F\left( a\right) =1$.  Thus
(\ref{KJ1}) follows from the previous part.
\end{proof}

\begin{lemma}
\label{ga} Suppose that $A$ is a uniform algebra on $\Omega,$ that $x_{0}$
is in the Choquet boundary of $A,$ and that $f$ and $p$ are strictly
positive continuous functions on $\Omega$. Then there is $a\in A$ such that $%
(af)\left( x_{0}\right) =p\left( x_{0}\right) $ and $\left| 
af \right| \leq p.$

Moreover only points from the Choquet boundary have this property.
\end{lemma}

\begin{proof}
Follows from \cite{Gam} II.12.
\end{proof}

\begin{lemma}
\label{Choquet} Suppose that $A$ is a uniform algebra on $\Omega$ and that $%
f $ is a strictly positive continuous functions on $\Omega$. Then $Ch\left(
A\right) \subset Ch(Af)\subset\Omega.$
\end{lemma}

\begin{proof}
To prove the first inclusion
assume that $x\in Ch\left( A\right)$, and assume
that $x=\frac{1}{2}F_{1}+\frac{1}{2}F_{2}$ where $F_{1},F_{2}$ are 
norm one functionals on $Af.$ Let $\mu _{1},\mu _{2}$ be norm one extensions
of $F_{1},F_{2}$ to functionals on $C\left( \Omega \right) $. By Lemma 
\ref{ga}, there is a net $a_{\gamma }$ in $A$ such that $\left\| a_{\gamma
}\right\| =a_{\gamma }\left( x\right) =1$ and $a_{\gamma }\rightarrow 0$
uniformly on compact subsets of $\Omega \backslash \left\{ x\right\} $. We
have
\begin{equation*}
f\left( x\right) = (a_{\gamma }f)\left( x\right) =\frac{1}{2}\int_{\Omega
}a_{\gamma }fd\mu _{1}+\frac{1}{2}\int_{\Omega }a_{\gamma }fd\mu
_{2}\rightarrow \frac{f\left( x\right) }{2}\left( \mu _{1}+\mu _{2}\right)
\left( \left\{ x\right\} \right) ,
\end{equation*}
hence $\mu _{1}\left( \left\{ x\right\} \right) =1=\mu _{2}\left( \left\{
x\right\} \right) $.
Since $\left\| \mu _{i}\right\| \leq 1$ we get $\mu
_{1}=x=\mu _{2}$.
\end{proof}

\vspace{2 mm}

Notice that, in general, $Ch(A)$ may be a proper subset of $Ch(Af)$. If, for
example, $A$ is equal to the disk algebra,
if $\Omega$ is the closed unit disk $\bar{\Di}$, and if
$f\left( z\right) = 2-\left| z\right| $, then $%
Ch\left( A\right) =\partial A=\partial\mathbb{D}$, while $0\in Ch(Af)$.

\begin{theorem}
\label{notAA} Assume that
$A$ and $B$ are uniform algebras on  compact sets $\Omega_1, 
\Omega_2$ respectively,
 and that
$f_{1},f_{2}$ are strictly positive continuous functions on $\Omega_1, 
\Omega_2$ respectively.   Suppose that
there is a surjective linear isometry $T:Af_{1}\rightarrow Bf_{2}$.
Then there is an invertible element $h$ of $B,$ and 
a homeomorphism $\varphi $
of $\Omega_2$ onto a subset $\varphi \left( \Omega_2 \right)$ of $M_{A}$,
 such that 
\begin{itemize}
\item [(i)]  the map $a \mapsto a \circ \varphi$  is an isometric 
isomorphism of $A$ onto $B$,
\item [(ii)]  $\varphi \left( \partial B\right) =\partial A$, 
\item [(iii)]  $\frac{f_{1}\circ \varphi }{f_{2}}_{|\partial
(Bf_{2})}= \left| h\right| _{|\partial (Bf_{2})}$, and
\end{itemize}

 \begin{equation}
T\left( af_{1}\right) =(a\circ \varphi) hf_{2}
\; \; \; \text{on} \; \Omega_2 \; , \; \text{for } \;
 a\in A .  \label{th1}
\end{equation}
Moreover, if $A = B$, then
 $T$ is an $A$-module isometry if and only if $\varphi $ is equal to
the identity map.
\end{theorem}

Before we prove this theorem, we give some
consequences.

\begin{corollary}
Assume that
$A$ is a uniform algebra on a compact set $\Omega$, and
that  $f$ is a
strictly positive continuous function on $\Omega$. Then $A$ and $Af$ are
linearly isometric if and only if they are $A$-isometric.
\end{corollary}

\begin{proof}
Assume that
$T:A\rightarrow Af$ is a linear isometry.  Set $f_{1}=1 , f_{2}=f$ 
in the last Theorem, and let $h$ be as in that Theorem.
 We can define a module map $%
S:A\rightarrow Af$ by
\begin{equation*}
S\left( a\right) \left( x\right) =ahf\left( x\right) \text{, }\qquad\text{%
for }a\in A,\quad x\in\Omega.
\end{equation*}
That $S$ is an isometry follows from the last Theorem and Lemma
\ref{ext} (iii).  \end{proof}

\vspace{2 mm}

In general it is not true that modules $Af_{1}$ and $Af_{2}$ are linearly
isometric if and only if they are $A$-module isometric.  Even if 
$f_1,f_2 \in C(\partial A)^+$, this is not true.
Assume, for example,
that $A$ is equal to the product of the disk algebra and the
two dimensional
algebra $C\left( \left\{ -1,1\right\} \right)$.  Set
$\Omega=\partial A=\mathbb{T%
}\times\left\{ -1,1\right\} $, and define a map
$\varphi : \Omega \rightarrow \Omega$ by
$\varphi\left( z,j\right) =\left( z,-j\right)
$, let $f\in C\left( \mathbb{T}\right) ^{+}\backslash Q(A(\Di))$, and put
\begin{equation*}
f_{1}\left( z,j\right) =\left\{
\begin{array}{ccc}
1 & \text{for} & j=-1, \\
f\left( z\right) & \text{for} & j=1,
\end{array}
\right.
\end{equation*}
and $f_{2}\left( z,j\right) =f_{1}\left( z,-j\right). $ 
Define $T\left(
af_{1}\right) \left( z,j\right)  = a\left( z,-j\right)
f_{2}\left( z,j\right) $.  Then $T$
 is a linear isometry from $Af_{1}$ onto $Af_{2}$.
Assume now that there is an $A$-module isometry $S:Af_{1}\rightarrow Af_{2}$. 
By Theorem \ref{notAA} there is an invertible element $h$ of $A$ such that
\begin{equation*}
S\left( af_{1}\right) =ahf_{2},\quad\text{for }a\in A.
\end{equation*}
Fix $\left( z,j\right) \in\mathbb{T}\times\left\{ -1,1\right\} $, and let $%
a_{n}$ be a sequence of norm one elements of $A$ convergent to zero almost
uniformly on $\mathbb{T}\times\left\{ -1,1\right\} \backslash\left\{
(z,j)\right\} $. The norm of $af_{1}$ is convergent to $\left| f_{1}\left(
z,j\right) \right| $, while the norm of $ahf_{2}$ is convergent to $\left|
h\left( z,j\right) f_{2}\left( z,j\right) \right| $ hence
\begin{equation*}
\left| h\left( z,j\right) f_{2}\left( z,j\right) \right| =\left| f_{1}\left(
z,j\right) \right| \text{, for }\left( z,j\right) \in\mathbb{T}\times\left\{
-1,1\right\} \text{.}
\end{equation*}
Hence $\left| f\left( z\right) \right| =\left| f_{1}\left( z,1\right)
\right| =\left| h\left( z,1\right) f_{2}\left( z,1\right) \right| =\left|
h\left( z,1\right) \right| $.  However $h\left( \cdot,1\right) $ is an
invertible element of the disk algebra,
 contrary to our assumption about $f$.

\begin{corollary}
\label{NEED}
Assume that
$A$ is a uniform algebra on a compact set $\Omega$, and that $f$ \ is 
a strictly positive continuous function on $\Omega$. Then the following
conditions are equivalent:

\begin{enumerate}
\item  $f\in Q$,

\item  $A,Af$, and $Af^{-1}$ are $A$-module isometric,

\item  $A,Af$, and $Af^{-1}$ are linearly isometric.
\end{enumerate}
\end{corollary}

\begin{proof}  The only thing that is still not clear here is that
($2 \Rightarrow 1$).  Note that if $T_1 : A \rightarrow Af$ is an
$A$-isometric isomorphism, and if
$T_1(1) = hf$, then $h \in A^{-1}$, $\Vert hf\Vert = 1$,
and $|h|f \leq 1$ on $\Omega$.
Since $\Vert a \Vert = \Vert ahf \Vert$ for all $a \in A$,
it follows by a simple Choquet point argument,
that $|h|f = 1$ on $Ch(A)$.   (Since if $x \in Ch(A)$
with $|h(x)|f(x) < \alpha < 1$, let $V = \{w \in \Omega :
|h(w)|f(w) > \alpha \}$.  By Lemma \ref{ga}, there exists $a \in A_1,
a(x) > \alpha,$ and $|a| < \frac{\alpha}{\Vert hf \Vert}$ on $V$.  Hence 
$\alpha < \Vert a \Vert = \Vert ahf \Vert \leq \alpha$.)

A similar argument shows that there exists $k \in A^{-1}$ with 
$|k| f^{-1} \leq 1$ on $\Omega$, and $|k| f^{-1} = 1$ on $Ch(A)$.
Hence $|hk| = 1$ on $Ch(A)$, and consequently on all of $\Omega$.
Thus $$f \leq |h^{-1}| = |k| \leq f \; $$
everywhere, so that $f = |k| \in Q$.
\end{proof}

\vspace{3 mm}

\begin{corollary}
\label{wasla}  Assume that
 $\Omega$ is a Shilov boundary of a uniform algebra $A$, and that $f$ is a 
strictly positive continuous functions on $\Omega$. Then the following
conditions are equivalent:

\begin{enumerate}
\item  $f\in Q$,

\item  $A$ and $Af$ are $A$-module isometric,

\item  $A$ and $Af$ are linearly isometric.
\end{enumerate}
Indeed, if $f_1, f_2 \in C(\partial A)^+$, then $Af_1 \cong Af_2$ 
$A$-isometrically, if and only if $f_1 f_2^{-1} \in Q$.
\end{corollary}

\begin{proof}  By Lemma \ref{Choquet}, $\Omega = \partial A = \partial(Af)$.
Now it is clear that Theorem \ref{notAA} gives ($3 \Rightarrow 1$).  
A similar argument proves the last statement.
The rest of the numbered equivalences are clear. 
  \end{proof}

\vspace{3 mm}

The last Corollary is not valid without the assumption that $\Omega$ is
equal to the Shilov boundary of $A.$ Indeed if we put $A=A\left( \mathbb{D}%
\right) $, $\Omega=\bar{\Di}$, and if we define
$f\left( z\right) = \frac{%
1+\left| z\right| }{2}$, then the map of multiplication by $f$ is an $A$%
-module isometry from $A$ onto $Af$.  On the other hand,
$f$ is not equal to an absolute
value of a function from $A$.

The first part of the last corollary also follows from:

\begin{corollary}
\label{new2}  Suppose that $A$ is a uniform algebra on $\Omega$, 
and $f \in C(\Omega)^+$.
Then $Af \cong A$ $A$-isometrically if and only if
 $Af \cong (Af)_{|(\partial A)}$ isometrically via the 
restriction map, 
and $f_{|(\partial A)} \in Q(\partial A)$.    
\end{corollary}

\begin{proof}  The $(\Leftarrow)$ direction is trivial.
The $(\Rightarrow)$ direction follows easily from the Theorem,
but we will give a different proof, which will generalize later 
to the `almost isometric' case.  
  Assuming that $Af \cong A$ $A$-isometrically,
as in the proof of Corollary \ref{NEED}, it 
follows that there exists a $h \in A^{-1}$ such that 
$\Vert a h f \Vert = \Vert a \Vert$, and that the latter 
statement implies that $|(hf)_{|(\partial A)}| = 1$ on 
$\partial A$, and $|hf| \leq 1$ on $\Omega$.  
Moreover for $w \in \Omega$, we have
$$|(af)(w)| \leq |(a h^{-1})(w)| \leq 
\Vert a h^{-1} \Vert_{\partial A} \leq
\Vert (af)_{|(\partial A)} \Vert \; . $$
  \end{proof}

\vspace{3 mm}

\begin{proof}[Proof of the Theorem.]
We will assume that $A = B$ for simplicity, although the same 
argument works in general.  Assume that
$T:Af_{1}\rightarrow Af_{2}$ is a linear surjective isometry. Let $%
K_{i}$, for $i=1,2$,
 be the set of extreme points of unit ball in the dual space $%
\left( Af_{i}\right) ^{\ast}$. Let $T_{|K_{2}}^{\ast}$ be a restriction of
the dual map $T^{\ast}$ to $K_{2}$. Since $T^{\ast}$ is a homeomorphism, in
the weak* topology, as well as being
an isometry, $T_{|K_{2}}^{\ast}$ is a
homeomorphism of $K_{2}$ onto $K_{1}$. Hence, for $\left| \alpha\right| =1$,
and $x\in Ch\left( Af_{2}\right) $ we have
\begin{equation*}
T^{\ast}\left( \alpha x\right) =\chi\left( \alpha,x\right) \cdot
\varphi\left( \alpha,x\right) \text{, where }\left| \chi\left(
\alpha,x\right) \right| =1\text{, and }\varphi\left( \alpha,x\right) \in
Ch\left( Af_{1}\right) \text{.}
\end{equation*}
Since $T^{\ast}$ is linear, the functions $\varphi$ and $\chi$ depend only
on $x,$ and we get the following representation of $T$:
\begin{equation*}
T\left( af_{1}\right) \left( x\right) =\chi\left( x\right) \cdot\left(
af_{1}\right) (\varphi\left( x\right)) \text{, \qquad for }a\in A,x\in
Ch\left( Af_{2}\right) ,
\end{equation*}
where $\varphi:Ch\left( Af_{2}\right) \rightarrow Ch(Af_{1})$ is a
surjective homeomorphism and $\chi$ is a unimodular function.

It is easy to see that this implies that $\varphi$ is continuous on
$Ch(Af_2)$.  (Since if $x_i \rightarrow x_0$ but $\varphi(x_i) \rightarrow
x_1 \neq \varphi(x_0)$, then choose $a \in A$ with $a(x_1) = 0 \neq
a(\varphi(x_0))$.  Then $T(af_1)(x_i) \rightarrow T(af_1)(x_0)$.
However $|T(af_1)(x_i)| \leq D |a(\varphi(x_i))|  \rightarrow 0$,
but $T(af_1)(x_0) \neq 0$.)

For $a=1$ we get that $Tf_{1}\left( x\right) =\chi\left( x\right) \cdot
f_{1}(\varphi\left( x\right)) =f_{2}(x) \frac{\chi\left( x\right) \cdot
f_{1}(\varphi\left( x\right)) }{f_{2}(x)}$.  Thus the
restriction of $\frac{\chi \cdot 
(f_{1}\circ\varphi)}{f_{2}}$ to $Ch(Af_2)$ is the restriction 
of a function 
in $A$.  Also  $\chi=\frac{%
Tf_{1}}{f_{1}\circ\varphi }.$ 
Hence it follows that $\chi$ is also continuous on $Ch(Af_2)$.

 Since $T $ is
surjective there is $a_{0}\in A$ such that $T\left( a_{0}f_{1}\right)
=f_{2}\cdot\left( \frac{\chi \cdot (f_{1}\circ\varphi)
}{f_{2}}\right) ^{2}$, as functions on $Ch\left( Af_{2}\right) .$ Hence
\begin{equation*}
f_{2}\cdot\left( \frac{\chi\cdot 
(f_{1}\circ\varphi)}{f_{2}}\right)^{2}=
T\left( a_{0}f_{1}\right) =f_{2}\cdot\frac{\chi\cdot(\left(
a_{0}f_{1}\right) \circ\varphi)}{f_{2}},
\end{equation*}
as functions on $Ch\left( Af_{2}\right) $, so
\begin{equation*}
\frac{\chi}{f_{2}}= (\frac{a_{0}}{f_{1}})\circ\varphi.
\end{equation*}
Consequently we get
\begin{equation}
T\left( af_{1}\right) \left( x\right) =f_{2}\left( x\right) \cdot (\left(
a_{0}a\right)( \varphi\left( x\right))) 
\text{, \qquad for }x\in Ch\left(
Af_{2}\right) \text{, }a\in A.  \label{brzeg}
\end{equation}
Let $b_{0}\in A$ be such that $T\left( b_{0}f_{1}\right) =f_{2}$. We have $%
f_{2}=T\left( b_{0}f_{1}\right) =f_{2}\cdot (\left( a_{0}b_{0}\right)
\circ\varphi)$, so $a_{0}b_{0}=1$ on $Ch\left( Af_{1}\right) $. By Lemma \ref
{Choquet} the closure of $Ch\left( Af_{1}\right) $ contains the Shilov
boundary of $A$ so $a_{0}b_{0}=1$ on $\partial A$, and consequently on $%
\Omega_1$.  This 
proves that $a_{0}$ is an invertible element of $A$. Hence $%
\left\{ a_{0}\cdot a:a\in A\right\} =A$.

Note that the map $a \longmapsto a\circ\varphi$ is linear, multiplicative,
one-to-one, and indeed is isometric since $Ch(A) \subset Ch(Af_1)$.   
As functions on $Ch(Af_2)$ we have:
\begin{equation}
\left\{ a\circ\varphi:a\in A\right\} =\left\{ \left( a_{0}\cdot a\right)
\circ\varphi:a\in A\right\} =T\left( Af_{1}\right) /f_{2}=A.  \label{brzeg2}
\end{equation}
Thus $a\longmapsto a\circ\varphi$ may be viewed as
an isometric automorphism $A \rightarrow A$.
In particular $a_{0}\circ\varphi$ is an invertible element of $A$ and we
have
\begin{equation*}
\frac{f_{1}\circ\varphi}{f_{2}}=\left| \frac{a_{0}\circ\varphi}{\chi}\right|
=\left| a_{0}\circ\varphi\right| \in Q
\end{equation*}
as required.

To finish the proof we need to show that $\varphi$ can be extended to a
homeomorphism of $\Omega_2$ onto a subset of $M_A$, and
that the formula (\ref{brzeg}) remains valid on the entire set $\Omega_2$. 
Since 
$a\longmapsto a\circ\varphi$ is an automorphism of $A$,
it is given by a
homeomorphism of the maximal ideal space of the algebra.  That is,
$\varphi$ can
be extended to a homeomorphism of $M_A$ onto itself
mapping $\partial A$ onto $\partial A$.   For a given $a$ in $A$, $%
T\left( af_{1}\right) $ and $f_{2}\cdot(\left( a_{0}a\right) \circ\varphi)$
are elements of $Af_{2}$ which, by
(\ref{brzeg}), are identical on the
Choquet boundary of $A$.  By dividing by $f_2$ if necessary,
we see that these elements are identical at any point
of $\Omega_2$ and we get (\ref{th1}).
\end{proof}

\vspace{2 mm}

Another proof of this result, using the function multiplier algebra
may be found in \cite{Bsh}.  However it is the proof above which extends 
to the `almost isometric' case.

\section{Almost isometries.}

Recall that Banach spaces $X,Y$ are almost isometric if the Banach-Mazur
distance between $X$ and $Y$ defined as
\begin{equation*}
d_{B-M}\left( X,Y\right) =\log \inf \left\{ \left\| T\right\| \left\|
T^{-1}\right\| :T:X\rightarrow Y\right\} ,
\end{equation*}
is equal to zero; two $A$-modules $X,Y$ are almost $A$-isometric if
\begin{equation*}
\log \inf \left\{ \left\| T\right\| \left\| T^{-1}\right\| :T:X\rightarrow
Y\right\} =0,
\end{equation*}
where this time the infimum is taken over the set of all $A$-isomorphisms.
Of course two isometric Banach spaces are almost isometric, but even for the
class of separable uniform algebras defined on subsets of a plane, almost
isometric spaces need not be isometric \cite{J1}. Small bound isomorphisms
between various classes of Banach spaces, primarily function spaces, have
been investigated in a large number of papers, see for example 
\cite{BC,J2,J3,RR}.

\begin{theorem}
\label{into}Assume 
that $A$ is a uniform algebra on a compact set $\Omega$, and that $f$
is a strictly positive continuous function on $\Omega$. Suppose that
$0 < \epsilon < \frac{1}{3}$ is given,
and suppose that $T$ is a
surjective linear map from 
$A$ onto $Af$ such that $\left\| T\right\| \leq 1$, 
and $\left\| T^{-1}\right\| \leq1+\varepsilon.$ 
 Then there is a subset $\Omega_{0}$ of
the Shilov boundary of $Af$ and a surjective continuous map $\varphi:\Omega
_{0}\rightarrow Ch(A)$ such that
\begin{equation}
\left| T\left( a\right) \left( x\right) -T\left( \mathbf{1}\right) \left(
x\right) \cdot a\circ\varphi\left( x\right) \right| \leq 4\varepsilon\frac{%
1+\varepsilon}{1-\varepsilon}\left\| a\right\| \text{, for }a\in A,\quad
x\in\Omega_{0} ,  \label{6.1}
\end{equation}

and such that 
\begin{equation}
1-10\varepsilon\leq\left| T\left( \mathbf{1}\right) \left( x\right) \right|
\leq1\text{, \quad for }x\in\Omega_{0} . \label{6.10}
\end{equation}
It follows that 
\begin{equation}
\sup\left\{ \left| (af)\left( x\right) \right| :x\in\Omega_{0}\right\}
\geq\left( 1- 15\varepsilon\right)
 \left\| af\right\| \text{,\quad\ for }af\in
Af,  \label{6.11}
\end{equation}
and it also follows that
the closure of $\Omega_{0}$ contains the Shilov boundary of $A$.
\end{theorem}

\begin{proof}
Let $\widetilde{T}:A\rightarrow Af_{|\partial(Af)}\subset C\left(
 \partial(Af)\right) $ 
be defined by $\widetilde{T}\left( a\right) = (1+\epsilon)
T\left( a\right)
_{|\partial(Af)}$.  We know from Lemma \ref{ext} (iii), that 
the restriction map from $Af$ to $Af_{|\partial(Af)}$ is an
isometry.  By Theorem 6.1 of \cite{J2} applied to $\widetilde{T}$,
 there is a subset $\Omega_{0}$
of $\partial(Af)$, and a continuous function $\varphi$ from $\Omega_{0}$ 
onto $Ch(A)$ such that (\ref{6.1}) holds .

Assume now that
there is an $x_{0}\in\Omega_{0}$ with $\left| T\left( \mathbf{1}%
\right) \left( x_{0}\right) \right| <1-10\varepsilon$. At the beginning of
the proof of Theorem 6.1 in \cite{J2}, the set $\Omega_{0}$ is defined
specifically as a subset of $\left\{ x\in\Omega:\left\| x\right\|
_{Af}> M\right\} $, where we denote by $\left\| x\right\| _{Af}$
the norm of the ``evaluation at the point $x$'' functional on $Af$,
and where $M$ can be chosen to
be any number satisfying $\frac{1-\varepsilon}{1+\varepsilon}>M>\frac
{1-\varepsilon}{1+\varepsilon}-\varepsilon^{2}$. Since $1-3\varepsilon
<\frac{1-\varepsilon}{1+\varepsilon}-\varepsilon^{2}$,  
there is a norm one element $a_{0}f$ of $Af$ such that $(a_{0}f)\left(
x_{0}\right) \geq1-3\varepsilon$. Put $a_{1}=T^{-1}\left( a_{0}f\right) $.
We have $\left\| a_{1}\right\| \leq\left\| T^{-1}\right\| \leq 1+\varepsilon$
and $T\left( a_{1}\right) \left( x_{0}\right) \geq1-3\varepsilon$,  while
\begin{equation*}
\left| T\left( \mathbf{1}\right) \left( x_0\right) \cdot a_{1}\circ
\varphi\left( x_{0}\right) \right| \leq\left| T\left( \mathbf{1}\right)
\left( x_0\right) \right| \left\| a_{1}\right\| \leq\left( 1-10\varepsilon
\right) \left( 1+\varepsilon\right) \text{.}
\end{equation*}
This contradicts (\ref{6.1}) and shows (\ref{6.10}).

To prove (\ref{6.11}) assume that
there is a norm one element $af$ of $Af$ such
that
\begin{equation*}
\sup \left\{ \left| af\left( x\right) \right| :x\in \Omega _{0}\right\}
<\left( 1-15\varepsilon \right) .
\end{equation*}
Put $b=T^{-1}\left( af\right) $, let $\widetilde{x}\in Ch(A)$ be such that $%
\left| b\left( \widetilde{x}\right) \right| =\left\| b\right\| \geq 1$, and
let $x_{1}\in \Omega _{0}$ be such that $\varphi \left( x_{1}\right) =%
\widetilde{x}$. By (\ref{6.10}) we have that
 \begin{equation*}
\left| T\left( \mathbf{1}\right) \left( x_1 \right) \cdot b\circ
\varphi \left( x_{1}\right) \right| \geq 1-10\varepsilon ,\text{ while }%
\left| Tb\left( x_{1}\right) \right| =\left| af\left( x_{1}\right) \right|
\leq 1- 15\varepsilon ,
\end{equation*}
which contradicts (\ref{6.1}) and shows (\ref{6.11}).

To finish the proof we need to show that (\ref{6.11}) implies $\partial
A\subset\overline{\Omega}_{0}$. To this end, 
choose a point  
$x_{0}\in Ch(A)$ which is not in
$\overline{\Omega}_{0}$.
W.l.o.g. we may assume that $\Vert f \Vert \leq 1$.   By
Lemma \ref{ga} there is an $a\in A$ such that
\begin{equation*}
\left\| a\right\| =1=a\left( x_{0}\right) \text{ and }\left| a\left(
x\right) \right| <\frac{1}{2}\min f\text{, for }x\in\overline{\Omega}_{0}.
\end{equation*}
By (\ref{6.11}) we have
\begin{align*}
\min f & \leq\left| (af)\left( x_{0}\right) \right| \leq\left\| 
af\right\| \\
& \leq\frac{1}{1-15\varepsilon}\sup\left\{ \left| (af)
\left( x\right) \right|
:x\in\Omega_{0}\right\} \\
& \leq\frac{1}{1-15\varepsilon}\frac{1}{2}\min f .
 \end{align*}
This contradiction proves that $\partial A\subset\overline{\Omega}_{0}$.
\end{proof}

\begin{corollary}
\label{got}
Assume that 
$A$ is a uniform algebra on a compact set $\Omega$, and that $f$ is a
strictly positive continuous function on $\Omega$. The following conditions
are equivalent:

\begin{enumerate}
\item[\textbf{1.}]  $A$ and $Af$ are almost $A$-isometric,

\item[\textbf{2.}]  $A$ and $Af$ are almost isometric.
\end{enumerate}

Moreover if $\Omega$ is equal to the Shilov boundary the above conditions
are also equivalent to:

\begin{enumerate}
\item[\textbf{3.}]  $f\in\overline{Q}^{+}.$
\end{enumerate}
\end{corollary}

\begin{proof}
That (3) implies (1) implies (2) is left to the reader.

Assume that $A$ and $Af$ are almost isometric, and let $T:A\rightarrow Af$
be such that $\left\| T\right\| \leq 1$ and $\left\| T^{-1}\right\| \leq
1+\varepsilon $. Since $T$ is surjective and $\frac{T\left( \mathbf{1}%
\right) }{f}\in A$, there is an $a_{0}\in A$ such that $Ta_{0}=\left( \frac{%
T\left( \mathbf{1}\right) }{f}\right) ^{2}f$.  By (\ref{6.1}) we get
\begin{align*}
\left\| \left( \frac{T\left( \mathbf{1}\right) }{f}\right) ^{2}f-T\left(
\mathbf{1}\right) \cdot a_{0}\circ \varphi \right\|_{\Omega_0}
 & \leq 4\varepsilon
\frac{1+\varepsilon }{1-\varepsilon }\left\| a_{0}\right\| \\
& \leq 4\varepsilon \frac{1+\varepsilon }{\left( 1-\varepsilon \right) ^{2}}%
\left\| \left( \frac{T\left( \mathbf{1}\right) }{f}\right) ^{2}f\right\| \\
& =4\varepsilon \frac{1+\varepsilon }{\left( 1-\varepsilon \right) ^{2}}%
\left\| \left( T\left( \mathbf{1}\right) \right) ^{2}f^{-1}\right\| \\
& \leq 4\varepsilon \frac{1+\varepsilon }{\left( 1-\varepsilon \right) ^{2}}%
\left\| f^{-1}\right\| ,
\end{align*}
so by (\ref{6.10}) we have
\begin{equation}
\left| T\left( \mathbf{1}\right) \left( x\right) -f\left( x\right) \cdot
a_{0}(\varphi \left( x\right)) \right| \leq 4\varepsilon \frac{%
1+\varepsilon }{\left( 1-10\varepsilon \right) \left( 1-\varepsilon \right)
^{2}}\left\| f^{-1}\right\| \left\| f\right\| \text{, for }x\in \Omega _{0}.
\label{KJ5}
\end{equation}
From (\ref{6.1}) and (\ref{KJ5}) we obtain
\begin{equation}
\left| T\left( a\right) \left( x\right) -f\left( x\right) \cdot \left(
aa_{0}\right)(\varphi \left( x\right)) \right| \leq \varepsilon
^{\prime }\left\| a\right\| \text{, for }a\in A,\quad x\in \Omega _{0},
\label{KJ6}
\end{equation}
where $\varepsilon ^{\prime }=4\varepsilon \frac{1+\varepsilon }{\left(
1-10\varepsilon \right) \left( 1-\varepsilon \right) ^{2}}\left\|
f^{-1}\right\| \left\| f\right\| +4\varepsilon \frac{1+\varepsilon }{%
1-\varepsilon }$.

Let $b_{0}\in A$ be such that $Tb_{0}=f$. If $\varepsilon $ is sufficiently
small, then
(\ref{KJ6}) gives 
\begin{align*}
\left\| \mathbf{1}-\left( a_{0}b_{0}\right) \circ \varphi \right\| _{\Omega
_{0}}& \leq \left\| f^{-1}\right\| \left\| f-f\cdot (\left( a_{0}b_{0}\right)
\circ \varphi) \right\| _{\Omega _{0}} \\
& =\left\| f^{-1}\right\| \left\| Tb_{0}-f\cdot (\left( a_{0}b_{0}\right)
\circ \varphi) \right\| _{\Omega _{0}} \\
& <\frac{1}{2}\; \; .
 \end{align*}

Hence $Re \left( a_{0}b_{0}\right) \left( x\right) \geq
\frac {1}{2}$ for
any $x\in Ch(A)$, and consequently for any $x$ in the maximal ideal space of $%
A $. It follows that $a_{0}$ and $b_{0}$ are invertible in $A$. \ By (\ref
{KJ5}) and (\ref{KJ6}) the function $\frac{T\left( \mathbf{1}\right) }{f}%
\cdot\frac{T\left( a_{0}^{-2}\right) }{f}$ is approximately equal, on $%
\Omega_{0} $, to $(a_{0}\circ\varphi)\cdot (a_{0}^{-1}\circ\varphi)=1$. 
Since $\bar{\Omega_0}$ contains $\partial A$,
the function $\frac{T\left( \mathbf{1}\right) }{f}\cdot\frac{T\left(
a_{0}^{-2}\right) }{f}-\mathbf{1}$ is approximately equal to zero on the
maximal ideal space of $A$.  Thus, if 
$\varepsilon$ is sufficiently small, $%
\frac{T\left( \mathbf{1}\right) }{f}\cdot\frac{T\left( a_{0}^{-2}\right) }{f}
$ is an invertible element of $A$.  Consequently $\frac{T\left( \mathbf{1}%
\right) }{f}$ is invertible. 

We can now define an $A$-module isomorphism $%
S:A\rightarrow Af$ by
\begin{equation*}
Sa=T\left( \mathbf{1}\right) a\text{, for }a\in A\text{.}
\end{equation*}
Fix $a\in A$ and let $\widetilde{a}\in A$ be such that $T\widetilde
{a}=Sa=T\left( \mathbf{1}\right) a$. By Theorem \ref{into}, $T\widetilde{a}$
is close, on $\Omega_{0}$, to $T\left( \mathbf{1}\right) \widetilde{a}%
\circ\varphi$.  It follows that $a\approx\widetilde{a}\circ\varphi$ on $%
\Omega_{0}$.  
Consequently by (\ref{6.11}) and the fact that $T$ is an `almost isometry',
 we get $%
\left\| a\right\| \approx\left\| \widetilde{a}\circ\varphi\right\| =\left\|
\widetilde{a}\right\| \approx\left\| T\widetilde{a}\right\| =\left\|
Sa\right\|$.  Thus $S$ is also an `almost isometry'.

Now assume that $\Omega$ is equal to the Shilov boundary. By (\ref{6.10}) we 
have that
$ \left| T\left( \mathbf{1}\right) \right| \approx1$ on $\overline{\Omega }%
_{0}=\partial A.$  As we proved before, $\frac{T\left( \mathbf{1}\right) }{f}$
is an invertible element of $A$.  It follows that $f\in\overline {Q}^{+}.$
\end{proof}

We have the following complement to the previous corollary
(cf. Corollary \ref{new2}):

\begin{corollary}
\label{new}  Suppose that $A$ is a uniform algebra on $\Omega$, and that
$f \in C(\Omega)^+$.   Then $Af \cong A$ almost $A$-isometrically if and 
only if $Af \cong (Af)_{|(\partial A)}$ isometrically via the restriction
map,
and $f_{|(\partial A)} \in \bar{Q}^+(\partial A)$.
\end{corollary}

\begin{proof}   This proceeds almost identically to the 
proof of \ref{new2}.  As in the proof of Corollary \ref{NEED}, it 
follows that for all $\epsilon > 0$,
there exists a $h_\epsilon \in A^{-1}$ such that 
$$(1-\epsilon)\Vert a \Vert  \leq 
\Vert a h_\epsilon f \Vert \leq (1+\epsilon) \Vert a \Vert \; \; ,$$ 
and that the latter 
statement implies that $|(h_\epsilon f)_{|(\partial A)}| \approx  1$ on 
$\partial A$, and $|hf| \leq 1+ \epsilon$ on $\Omega$.  
Thus $f_{|(\partial A)} \in \bar{Q}^+(\partial A)$.
Moreover for $w \in \Omega$, we have
$$|(af)(w)| \leq (1+ \epsilon) |(a h_{\epsilon}^{-1})(w)|
 \leq (1+ \epsilon) \Vert a h_{\epsilon}^{-1} \Vert_{\partial A}
\leq (1 + \epsilon)^2 \Vert (af)_{|(\partial A)} \Vert \; 
. $$
Since $\epsilon > 0$ is arbitrary, $\Vert af \Vert = 
\Vert af \Vert_{(\partial A)}$.
  \end{proof}

\begin{corollary}
\label{notA} For any strictly positive $f\in C(\Omega )$, the following are
equivalent.
\begin{itemize}
\item[(i)]  $f\in \bar{Q}^{+}$,

\item[(ii)]  $A\cong Af\cong Af^{-1}$ almost $A$-isometrically,

\item[(iii)]  $A\cong Af\cong Af^{-1}$ linearly almost isometrically.
\end{itemize}
\end{corollary}

\begin{proof}
We will only  prove $(ii)\Rightarrow \left( i\right)$.  That $(iii)$ implies
$(ii)$ follows from what we just did.  The easy implications 
$(i) \Rightarrow (ii) \Rightarrow (iii)$ are left to the reader.

If $T_{1}:A\rightarrow Af$ and
$T_{2}:A\rightarrow Af^{-1}$ are module
isomorphisms, let $t_1 = T_{1}(\mathbf{1})$ and $t_2 =
T_{2}(\mathbf{1})$.
Then $T_{1}(a)=t_{1} a$, $%
T_{2}(a)=t_{2} a$, and $t_{1} f^{-1}$, $t_{2} f$ are 
invertible elements of $A$.
If $\left\| T_{1}\right\| \leq 1+\varepsilon ,\left\| T_{2}\right\| \leq
1+\varepsilon ,$ then
\begin{equation}
\left| t_{1}\left( x\right) \right| \leq
1+\varepsilon \quad \text{and\quad }\left| t_{2}
\left( x\right) \right| \leq 1+\varepsilon ,\quad \text{ for }x\in \Omega
\label{KJ7}
\end{equation}
and consequently
\begin{equation*}
\left| t_{1}\left( x\right) t_{2}\left( x\right) 
\right| \leq \left( 1+\varepsilon \right) ^{2},\quad
\text{ for }x\in \Omega .
\end{equation*}
Suppose that
 $\left\| T_{1}^{-1}\right\| \leq 1,\left\| T_{2}^{-1}\right\| \leq 1$.
If $\left| t_{1}\left( x_{0}\right) \right| =1-r<1$
 for some $x_{0}\in Ch(A)$, then
by Lemma \ref{ga} there
is $a\in A$ such that $a\left( x_{0}\right) =1=\left\| a\right\| $ and $\left|
a\left( x\right) \right| <\frac{1}{2+\varepsilon }$ for $x\in \left\{ x:\left|
t_{1}\left( x\right) \right| >1-\frac{r}{2}\right\} $. 
We get $\left\| T_{1}\left( a\right) \right\|
=\left\| t_{1} a\right\| \leq \max \left\{ \frac{1}{2+\varepsilon
}\left\| T_{1}\right\| ,1-\frac{r}{2}\right\} <1$.
This contradicts the assumption
that $\left\| T_{1}^{-1}\right\| \leq 1$.

Hence
\begin{equation*}
\left| t_{1}\left( x\right) \right| \geq 1\quad
\text{and\quad }\left| t_{2}\left( x\right) \right|
\geq 1,\quad \text{ for }x\in Ch(A),
\end{equation*}

so
\begin{equation*}
\left| t_{1} \left( x\right) t_{2}\left( x\right) 
\right| \geq 1,\quad \text{ for }x\in Ch(A).
\end{equation*}
However $t_{1} t_{2}$ is
invertible, so it attains minimum on $Ch(A)$.  Hence
\begin{equation*}
1\leq \left| t_{1} \left( x\right) t_{2}
\left( x\right) \right| \leq \left( 1+\varepsilon 
\right)^{2},\quad \text{ for }x\in \Omega .
\end{equation*}
By (\ref{KJ7}) it follows that
\begin{equation*}
\frac{1}{1+\varepsilon }\leq \left| t_{2} \left(
x\right) \right| \leq 1+\varepsilon ,\quad \text{ for }x\in \Omega .
 \end{equation*}
Thus $\left\| f-\left| t_{2} \right| f \right\| \leq
\varepsilon \left\| f\right\| $; and
since $\varepsilon >0$ is arbitrary we get $%
f\in \bar{Q}^{+}$.
\end{proof}

\begin{lemma}  \label{prev}  If $X$ is a Banach $A$-module over 
a function algebra $A$, and if $X$ 
is almost $A$-isometric to a function $A$-module, then 
$X$ is a function $A$-module.
\end{lemma}

\begin{proof}  We use the 
 injective tensor norm characterization of 
function modules  \cite{BLM}.
Suppose that $T_\epsilon$ is the $\epsilon$-isomorphism.  
Then for $a_1, \cdots , a_n \in A$, and
$x_1, \cdots , x_n \in X$, we have:
$$ \Vert \sum_i a_i T_\epsilon(x_i) \Vert \leq 
\Vert \sum_i a_i \otimes T_\epsilon(x_i) \Vert_\lambda
\leq \Vert T_\epsilon \Vert \Vert \sum_i a_i \otimes x_i\Vert_\lambda \; \; ,
$$
where $\lambda$ is the injective tensor norm.
Thus:
$$ \Vert \sum_i a_i x_i \Vert = \lim_{\epsilon \rightarrow 0}
\Vert T_\epsilon (\sum_i a_i x_i) \Vert 
\leq \Vert \sum_i a_i \otimes x_i\Vert_\lambda \; \; . $$
By \cite{BLM}, $X$ is a function $A$-module.
\end{proof}

\begin{corollary}
\label{Motat21}  Let $X$ be a Banach
$A$-module.  The following are equivalent:
\begin{itemize}
\item [(i)]  $X$ is an algebraically singly generated, faithful
function $A$-module, and $X \cong A$ almost linearly isometrically,
\item [(ii)]  $X \cong A$ almost $A$-isometrically,
\item [(iii)] There exists $f \in \bar{Q}^+$ such that $X \cong Af$
$A$-isometrically.
\end{itemize}
\end{corollary}

\begin{proof}  Assume (i).  As we said early in the introduction,
the first two conditions in (i), together with
a corollary in \S 3 of \cite{Bsh}, shows that $X \cong Af$
$A$-isometrically, where $f \in C(\Omega)^+$, for some 
$\Omega$ on which $A$ acts as a function algebra.  
Now (ii) follows from Corollary \ref{got}.  By Corollary 
 \ref{new},
$Af \cong (Af)_{|(\partial A)}$ and $f_{|(\partial A)} \in 
\bar{Q}^+(\partial A)$, showing (iii).
   
Given (ii), it follows by Lemma \ref{prev} that $X$ is a function
$A$-module, and now (i) is clear.   
Clearly  (iii) implies (i).
\end{proof}

\vspace{2 mm}

A similar result and proof holds with almost isometries replaced 
by isometries.
 
We noticed that the results in this section
 have been stated
for a single function $f$, rather then for a pair of functions $f_{1},f_{2}$,
like results in the previous section concerning isometries. 
There are analogous results describing almost
isometries between modules $Af_{1}$ and $Af_{2}$,
 however they are
much more technical and involve not a single automorphism $\varphi :\Omega
\rightarrow \Omega $ but a sequence of homeomorphisms between the Choquet
boundaries of $Af_{1}$ and $Af_{2}$. The following theorem, which may be of
an independent interest, can be proven using methods similar to that of the
proof of Theorem 6.1 of \cite{J2}.  It may then be used to extend the last
few results to almost isometries between modules $Af_{1}$ and $Af_{2}.$

\begin{theorem}
Assume that $\Omega$ is the Shilov boundary of a uniform algebra $A$, and 
that $%
f_{1},f_{2}$ are strictly positive continuous functions on $\Omega$. Suppose
that
 $T:Af_{1}\rightarrow Af_{2}$ is a surjective linear isomorphism such that $%
\left\| T\right\| \leq1+\varepsilon$ and $\left\| T^{-1}\right\|
\leq1+\varepsilon$ where $\varepsilon<\varepsilon_{0}$ (an absolute
constant). Then there is a dense subset $S$ of $\Omega$ containing the
Choquet boundary of $A,$ a continuous bijection $\varphi:S\rightarrow Ch(A),$
and a continuous unimodular function $\chi$ such that
\begin{equation}
\left\| T\left( af_{1}\right) -\chi\cdot\left( af_{1}\right) \circ
\varphi\right\| \leq\varepsilon^{\prime}\left\| af_{1}\right\| \text{,\quad\
for }a\in A\text{,}
\end{equation}
where $\varepsilon^{\prime}\rightarrow0$ as $\varepsilon\rightarrow0$.
\end{theorem}

If the functions $f_1,f_2 \in \M_A$, then a simple version of an
 `almost isometry' result may be found in Theorem \ref{gen}.

No doubt there are also versions of all these results for a pair 
of function algebras $A$ and $B$, and function modules 
$Af_1$ and $Bf_2$, but that will take us a little further 
afield from our main concerns.

\begin{center}
Part B.
\end{center}

\section{Some observations concerning approximations in modulus}

In this section we make various observations 
concerning the sets $P,Q,\mathcal{G},%
\mathcal{F,} G_{A}(\Omega), H_{A}(\Omega )$
and $\M_A(\Omega)$
 defined in the Introduction,  Notation
and Definitions section.
These sets will play
crucial role in the next section when we study singly generated bimodules.

We already observed that $C(\Omega )^{+}
= \bar{Q}^{+}$
if and only if $A$ is logmodular on $\Omega $. 
 The question of when $C(\Omega )_{+}=%
\bar{P}$ (resp. $C(\Omega )_{+}=\mathcal{G}$) has been studied by Mlak,
Glicksberg, Douglas and Paulsen, and others. An algebra with this property
is called `approximating in modulus' (resp. `convexly approximating in
modulus').   Just as in the usual proof 
that $A$ may only be logmodular on $\partial A$, one can show
 using Urysohn's lemma,
Lemma \ref{ga},
and the fact that a function in $A$ has maximum modulus achieved on
$\partial A$, that either of these `approximating in modulus'
properties forces
 $\Omega = \partial A$. 

For example, Dirichlet algebras, such as the disk algebra $A(%
\mathbb{D})$ considered as functions on the circle, are logmodular, and
approximating in modulus. Glicksberg gave the following sufficient
condition: If the inner functions (that is, functions in $A$ which have
constant modulus 1 on $\Omega $) separate points of $\Omega $, then $A$ is
approximating in modulus.

By the Stone-Weierstrass theorem it is easy to see that for any function
algebra $A$ on $\Omega $ the set $\mathcal{F}-\mathcal{F}$ is
dense in $%
C_{\Rl}(\Omega)$. However this does not imply that $\mathcal{F}$ is
dense in $C(\Omega )_{+}$, or equivalently, that $A$ is convexly
approximating in modulus. Indeed we have the following:

\begin{proposition}
Let $A$ be a function algebra. If $Ch(A)\neq\partial A$ then $A$ is not 
convexly approximating in modulus on $\partial A$.  That is,
$\mathcal{G}\neq
C\left( \partial A\right) _{+}$.
\end{proposition}

\begin{proof}
Let $x_{0}\in \partial A,$ let $V\subset \partial A$ be an open neighborhood
of $x_{0},$ and let $f\in C\left( \partial A\right) _{+}$ be such that
\begin{equation*}
f\left( x_{0}\right) =1=\left\| f\right\| \text{, and }f\left( x\right) =0%
\text{ for }x\in \partial A-V\text{.}
\end{equation*}
Assume that $\mathcal{G}=C\left( \partial A\right) _{+}$, and let $%
g_{1},...,g_{n}\in A$ be such that
\begin{equation*}
\left| f-\sum_{j=1}^{n}\left| g_{j}\right| ^{2}\right| \leq \frac{1}{4}.
\end{equation*}
Multiplying $g_{1},...,g_{n}$ by appropriate numbers of absolute value $1$
we may assume that $g_{j}\left( x_{0}\right) \in \mathbb{R}_{+}$.
Then \begin{equation*}
\frac{3}{4}\leq \left( \sum_{j=1}^{n}\left| g_{j}\right| ^{2}\right) \left(
x_{0}\right) =\sum_{j=1}^{n}g_{j}\left( x_{0}\right) ^{2}\leq \frac{5}{4}.
\end{equation*}
Put
\begin{equation*}
g=\frac{\sum_{j=1}^{n}g_{j}^{2}}{\sum_{j=1}^{n}g_{j}\left( x_{0}\right) ^{2}}%
\in A.
\end{equation*}
We have
\begin{equation*}
g\left( x_{0}\right) =1\text{, }\left\| g\right\| \leq \frac{5}{3}\text{,
and }\left| g\left( x\right) \right| \leq \frac{1}{3}\text{ for }x\in
\partial A-V\text{.}
\end{equation*}
By the Bishop ''$\frac{1}{4}-\frac{3}{4}$'' criterion (\cite{Gam} Th. 11.1
page 52 and remark on p. 59) $x_{0}$ is a p-point of $A$. 
That is $x_{0}\in Ch(A)$.
\end{proof}

\vspace{3 mm}

\begin{lemma}
\label{Gmax}  Suppose that $f^2 \in \G(\Omega)$.  Then $f$ 
achieves its norm on $\partial A$.  Indeed $Af \cong
(Af)_{|(\partial A)}$ $A$-isometrically.
\end{lemma}

\begin{proof}  Suppose that $K_n$ are $A$-tuples with $\Vert K_n(\cdot)
\Vert_2$ converging uniformly to $f$.  Given $\epsilon > 0$, 
we have $\Vert K_n(w) \Vert_2 \leq f(w) + \epsilon$ for all
$w \in \Omega$ and sufficiently large $n$.  For such $w , n$, and
for any complex Euclidean unit vector $z$ we have
$$|z.(a(w)K_n(w))| \leq \Vert z.(a(\cdot)K_n(\cdot)) \Vert_{\partial A}
   \leq \Vert af \Vert_{\partial A} + \epsilon \Vert a \Vert \; \; .
$$
Thus $\Vert a(w)K_n(w) \Vert_2 \leq 
\Vert af \Vert_{\partial A} + \epsilon \Vert a \Vert$.
Letting $n \rightarrow \infty$ gives
$\Vert af \Vert_{\Omega} \leq
\Vert af \Vert_{\partial A} + \epsilon \Vert a \Vert$.
Since $\epsilon > 0$ was arbitrary, we get the result.
 \end{proof}

Using this we can sharpen an earlier result:

\begin{corollary}
\label{well}  Suppose that $f \in C(\Omega)^+$,
that $f^{-2} \in \G(\Omega)$, and that $Af \cong A$ 
$A$-isometrically (resp. almost $A$-isometrically).  Then
$f \in Q(\Omega)$ (resp. $f \in  \bar{Q}^+(\Omega)$).
\end{corollary}

\begin{proof}  By the Lemma and Corollary \ref{new2} 
(resp. \ref{new}), $f_{|(\partial A)}
\in Q(\partial A)$ (resp. in $\bar{Q}^+(\partial A)$), and $Af^{-1} \cong 
(Af^{-1})_{|(\partial A)} \cong A$.   Now use
corollary \ref{NEED} (resp. \ref{notA}).   
\end{proof}

\vspace{3 mm}

Note that for the disk algebra $A(\mathbb{D})$ considered as functions on
the closed disk, or more generally for any function algebra containing no
nontrivial inner functions it is clear that $P\cap P^{-1}=Q$. Thus for
certain good examples, we will have to look among algebras with many inner
functions.   It is also easy to see that $P(M_A) \cap P^{-1}(M_A)
= Q(M_A)$.

It will be significant to us that $Q,\bar{Q}^{+}$ and $\M$
 are (abelian) groups.  We will write
$Q'$ for the quotient group $\bar{Q}^{+}/Q$.
Of course for some function algebras,
it can happen that $Q=\bar{Q}^{+}$,
for example for $H^{\infty }(%
\mathbb{D})$, where $\bar{Q}^{+} = L^{\infty }\left( \mathbb{T}\right)^{+}$
(\cite{Su}, Th. 5.26). On the other hand, for $A(\mathbb{D})$,
the disk algebra, we have that $Q$ is a proper subset of
$\bar{Q}^{+}=C(\mathbb{T})^{+}$.

A good way of producing interesting functions in $\bar{Q}^{+}$, which works
in any nonselfadjoint function algebra goes as follows: By a result of
Hoffman and Wermer \cite{sto}, $Re A$ is not uniformly closed. If $%
g\in (Re \; A)^{\bar{}}\setminus Re A$, set $f=e^{g}$. If $a_{n}\in
A,\; Re a_{n}\rightarrow g$ uniformly, then $e^{Re \; %
a_{n}}=|e^{a_{n}}|\rightarrow f$. Hence $f\in \bar{Q}^{+}$. If the set of
invertible elements in $A$ is connected, or equivalently (\cite{Gam} p. 91),
if the first Cech cohomology group $H^{1}(M_{A},\mathbb{Z})=0$ of the
maximal ideal space $M_{A}$ is zero, then $h\in Q$ if and only if $h=e^{%
Re \; a}$ for some $a\in A$; thus we can definitely assert that $f\notin
Q $ in this case.  Thus we have proved:

\begin{corollary}
\label{con} Suppose that $A$ is a nonselfadjoint function algebra, such that
$H^{1}(M_{A},\mathbb{Z})=0$, then $Q^{\prime }=\bar{Q}^{+}/Q\neq 0$. Hence
there exist nontrivial function $A$-modules which are almost
$A$-isometric to $A$.
\end{corollary}

We shall see later that this corollary also gives 
the existence of
nontrivial rank 1 strong Morita equivalence bimodules for any such
$A$.

We now turn to the class $\M_A(\Omega)$ defined in the introduction.  We have
  \begin{equation*}
G_{A}(\Omega )\subset H_{A}(\Omega )\subset \log \mathcal{M}_{A}(\Omega ).
\end{equation*}
In Part C we shall see that
all three of the above classes have the five properties of `harmonic
classes' with respect to $A$ described in the introduction.
We shall not use the following, but state it for interests sake.
Its proof follows from the definition of $\M$ in the introduction, and
is left to the reader.

\begin{proposition} 
If  $f \in \M_A(\Omega)$ 
and if $t > 0$ , then $tf \in \M$.  Moreover,
$\M$ is closed in the relative topology from $C(\Omega)^+$.
Thus $\log \M_A(\Omega)$ is uniformly closed.
\end{proposition} 

 Note that $\M$ is
not complete in the norm topology (since if $f \in \M$ and
$t > 0$ , then $tf \in \M$, but $\lim_{t \rightarrow 0} tf = 0
\notin \M$).

The authors do not know if $H_{A}(\Omega )=\log \mathcal{M}_{A}(\Omega )$ in
general. In any case, either answer to the question seems very interesting.
If $H(\Omega )=\log \mathcal{M}_{A}(\Omega )$ then we obtain from what 
we do later, amongst other
things, a neat description of all topologically singly generated Morita
equivalence $A-A$-bimodules; but if $H(\Omega )\neq \log \mathcal{M}%
_{A}(\Omega )$ in general, then $\log \mathcal{M}_{A}(\Omega )$ seems to be
a genuinely new and interesting class of harmonic functions with respect to $%
A$.   
A uniform algebra is called `Dirichlet' if 
$G_{A}(\Omega )=$ $C_{\Rl}(\Omega)$, is called 
`logmodular' if $H_{A}(\Omega )= C_{\Rl}(\Omega)$;
therefore
 we will call an algebra 
\emph{logMorita} if $\log \mathcal{M}_{A}(\Omega )=
C_{\Rl}(\Omega)$. Every logmodular algebra is logMorita.
The next result shows that logMorita algebras share many properties
with logmodular algebras.

\begin{theorem}
\label{Omes} If $A$ is a logMorita function algebra on $\Omega$, then

\begin{itemize}
\item[(i)]  $\Omega $ is the Shilov boundary of $A$,

\item[(ii)]  Every $\phi \in M_A$ has a unique representing measure,

\item[(iii)]  If $\Pi $ is a Gleason part for $A$, then either $\Pi $ is a
singleton, or is an analytic disk, in the sense that there is a bijective
continuous map $\Phi :\mathbb{D}\rightarrow \Pi $ , such that if $f\in A$
then $\hat{f}\circ \Phi $ is holomorphic on $\mathbb{D}$.
\end{itemize}
\end{theorem}

\begin{proof}
(i) is proved as for Dirichlet and logmodular algebras, using Urysohns lemma
and the fact that for an
 $A$-tuple $H \in A^{(n)}$, the function $\Vert H(\cdot )\Vert_{2}$ 
achieves its maximum modulus on the Shilov boundary.   (The latter fact
may be seen by considering $z.H(w)$ for $z \in \Co^n_1$).
Similarly, (ii)
follows the classical line of proof: Suppose that $\mu ,\nu $ are
representing measures for $\phi \in M_{A}$, and that $H,K$ are $A$-tuples
with $1=H(w)\cdot K(w)$ for all $w\in \Omega $. By Fubini and
Cauchy-Schwarz, we have
\begin{equation*}
1=\hat{H}(\phi )\cdot \hat{K}(\phi )=\int_{\Omega \times \Omega }H(w)\cdot
K(z)\;d(\mu \times \nu )\leq (\int \Vert H(w)\Vert _{2}d\mu )(\int \Vert
K(w)\Vert _{2}d\nu )\;.
\end{equation*}
The remainder of the proof follows 17.1 in \cite{sto}. Finally, (iii)
follows from (ii) by 17.1 of \cite{sto}.
\end{proof}

\section{Approximations in modulus and equivalence bimodules}

In \cite{BMP} the notion of \emph{strong Morita equivalence} is defined for
a pair of operator algebras $A$ and $B$. Its theory and consequences have
been worked out there and in other papers of ours (see \cite{BMN,BOMD} for
example). A related notion, \emph{strong subequivalence}, was recently
defined in \cite{BOMD}.  It was shown to have many of the properties 
of strong Morita equivalence.  One of our 
objectives here is to show that it is not
the same as strong Morita equivalence.
It will not
be necessary for us to state the general definitions of these notions here,
we will simply say that they involve a pair of bimodules $X$ and $Y$, called
\emph{equivalence bimodules}. In fact we
shall restrict our attention here to the special case where the operator
algebras are function algebras, and in this case the two definitions can be
simplified.  Indeed,  for a bimodule of the form 
$Af$, where $f$ is a strictly positive function
on $\Omega$, and considering the canonical pairing $Af \times Af^{-1} 
\rightarrow A$, it is easy to translate the definitions from \cite{BMP,BOMD},
using Lemma 2.8 of \cite{BMP}, into the following precise form:

\begin{definition}
\label{prec} Suppose that $A$ is a function algebra on a compact space $%
\Omega$, and $f$ is a strictly positive continuous function on $\Omega$.

\begin{itemize}
\item[(i)]  We say that 
$Af$ is a \emph{strong Morita equivalence bimodule
(with inverse bimodule $Af^{-1}$)},
if whenever $\epsilon >0$ is given, then we can write $1=%
\sum_{i=1}^{n}x_{i}y_{i}$ as functions on $\Omega $, with $x_{i}\in
Af,y_{i}\in Af^{-1}$, and $\sum_{i=1}^{n}|x_{i}(\omega )|^{2}\leq 1+\epsilon
$, and $\sum_{i=1}^{n}|y_{i}(\omega )|^{2}\leq 1+\epsilon $, for all $\omega
\in \Omega $. We will say $Af$ is rank 1, if $n=1$ in the above.

\item[(ii)]  We say that $Af$ is a \emph{strong subequivalence bimodule} if
whenever $\epsilon >0$ is given, there are $x_{i}\in Af,y_{i}\in Af^{-1}$,
such that $1-\epsilon \leq \sum_{i=1}^{n}|x_{i}(\omega )|^{2}\leq 1+\epsilon,
$ and $1-\epsilon \leq \sum_{i=1}^{n}|y_{i}(\omega )|^{2}\leq 1+\epsilon $,
for all $\omega \in \Omega $.

\item[(iii)]  We say that $Af$ is a \emph{unitary subequivalence bimodule}
if (ii) holds, but with $\epsilon =0$ and $n=1$.

\item[(iv)]  We shall say that $Af$ is a \emph{Shilov subequivalence 
bimodule} 
if it is a strong subequivalence bimodule and $\Omega$ is the Shilov
boundary of $A$.
\end{itemize}
\end{definition}

All these definitions are in \cite{BMP,BOMD} except (iii). Also, strictly
speaking, in (i) we should say that $(A,A,Af,Af^{-1},\cdot,\cdot)$ \emph{is
a strong Morita context} (see \cite{BMP} Definition 3.1). Here `$\cdot$'
refers to multiplication of scalar functions on $\Omega$. However, to avoid
this somewhat cumbersome notation, we will use the looser convention of (i).
Also, concerning (iv), we used the word `minimal' instead of `Shilov' in
\cite{BOMD} Definition 5.7.

It is clear in the definitions above,
that $(iii)\Rightarrow (ii)$ and that $(i)\Rightarrow (ii)$. 

\begin{proposition}
\label{ma} Suppose that $f\in C(\Omega )^{+}$. Then:

\begin{itemize}
\item[(a)]  $Af$ is a unitary subequivalence bimodule if and only if $f\in
P\cap P^{-1}$. That is, iff $f\in P$ and $f^{-1}\in P$.

\item[(b)]  $Af$ is a strong subequivalence bimodule if and only if $%
f^{2}\in \mathcal{G}\cap \mathcal{G}^{-1}$

\item[(c)]  $f\in \bar{Q}^{+}$ if and only if $Af$ is a 
rank 1 strong Morita equivalence bimodule (with inverse 
bimodule $Af^{-1}$).

\item[(d)]  $f \in \M(\Omega)$ if and only if $Af$ is a
strong Morita equivalence bimodule (with inverse 
bimodule $Af^{-1}$).
\end{itemize}
\end{proposition}

\begin{proof}
\begin{itemize}
\item[(a)]  Clearly $Af$ is a unitary subequivalence bimodule if
and only if $Af$ and $Af^{-1}$ both contain inner functions (functions of
constant modulus 1 on $\Omega $), which is clearly equivalent to $f\in P\cap
P^{-1}$.

\item[(b)]  Suppose that $\epsilon >0$ is given, and that $%
x_{i},y_{i}$ are as in (ii). Then $\Vert 1-\sum_{i}|x_{i}|^{2}\Vert _{\infty
}\rightarrow 0$ as one allows $\epsilon \rightarrow 0$. Hence $f^{-2}\in
\mathcal{G}$. Similarly $f^{2}\in \mathcal{G}$. For the converse, notice
that the argument given is reversible.

\item[(c)]  If $|g_{m}|\rightarrow f$ uniformly as $m\rightarrow
\infty $, where $g_{m},g_{m}^{-1}\in A$, then we may write $%
1=(g_{m}^{-1}f)(g_{m}f^{-1})$, and clearly $\Vert g_{m}^{-1}f\Vert _{\infty
}\rightarrow 1$ and $\Vert g_{m}f^{-1}\Vert _{\infty }\rightarrow 1$. This
obviously implies the condition in (i) of Definition \ref{prec} with $n=1$.
Also, this argument is reversible.

\item [(d)]  This is obvious.
\end{itemize}
\end{proof}

\vspace{4 mm}

\noindent {\bf Remark.}  To those readers who are familiar 
with notions in \cite{BOMD}, we remark that 
the proof of (b) shows the
following.  Consider $Z = C(\Omega)$ as a module over itself.
It is clearly generated as a $C(\Omega)$-module by the
$A$-submodule $Af$.  The proof of (b) shows,
in the language of \cite{BOMD}, that  $Z$
is the {\em C$^*-$dilation}
of $Af$ if and only if $f^{-2} \in \G$.

\begin{theorem}
\label{oop}
Suppose that $A$ is a function algebra which is antisymmetric on
$\Omega$, or for which the number of pieces in the antisymmetric
decomposition\footnote{See \cite{sto} Theorem 12.1 for example.
An algebra with finite antisymmetric decomposition, is a
direct sum of finitely many uniform algebras.}
 of $A$ is finite.  Suppose that
$f$ is a function such that $Af$ satisfies
Definition \ref{prec} (i) with $\epsilon = 0$.  Then $\log f \in
H(\Omega)$, and indeed $f \in Q$.
\end{theorem}

\begin{proof}
First suppose that $A$ is antisymmetric on
$\Omega$.   By hypothesis, we have $f = \Vert K(\cdot) \Vert_2$ and
$f^{-1} = \Vert H(\cdot) \Vert_2$, for $A$-tuples $H,K$ with
$1 = H.K = (Hf).(Kf^{-1})$.   By the converse to Cauchy-Schwarz,
$Hf = K^*f^{-1}$, where `$*$' is the complex conjugate.
Thus $h_ik_i = |k_i|^2 f^{-2} \geq 0$.
By antisymmetry, $h_i k_i = c_i$ a nonnegative constant.
Therefore, if $c_1 \neq 0$ then $h_1,k_1 \in A^{-1}$.
We have $c_1 = h_1 k_1 = |k_1|^2 f^{-2}$, so that
$f = \frac{1}{\sqrt{c_1}} |k_1| \in Q$.

If $\Omega$ is a disjoint union of a finite number of 
antisymmetric
pieces $\Omega_i$, then each $\Omega_i$ is open and compact,
 and by the first part we have
$f_{|\Omega_i} = |b_i|$, for an invertible
$b_i \in A_{|\Omega_i}$.   Put $b(x) = b_i(x)$ if $x \in
\Omega_i$, then $b \in C(\Omega)$.  It
 follows from \cite{sto} Theorem 12.1 say, that $b \in A$
and $b^{-1} \in A$.   Thus $f = |b| \in Q$.
\end{proof}

\vspace{3 mm}
 
Comparing \ref{ma} (c) and \ref{wasla} shows that for any 
$f\in \bar{Q}^{+}\setminus Q$, provided say by \ref{con},
we have that $Af$ is
 a nontrivial strong Morita equivalence bimodule. 

From the next theorem we will be able to
give examples of strong subequivalence bimodules which
are not strong Morita equivalence bimodules.

\begin{theorem}
Let $A$ be a uniform algebra on $\Omega$.  Suppose that
$w_{1},w_{2}$ are distinct
points of $\Omega$, and that $G$ is 
an invertible element of $A$ such that $G\left(
w_{1}\right) =1=-G\left( w_{2}\right) $. Put $A_{0}=\left\{ a\in A:a\left(
w_{1}\right) =a\left( w_{2}\right) \right\}$, and 
set $f=\left| G\right| $. Then
$f \in \M_{A_0}$ if and only if $w_{1},w_{2}$
are in different Gleason parts of $A$, and  if and only if
$f \in \bar{Q}^+_{A_0}$.
\end{theorem}

\begin{proof}
Assume that $f \in \M_{A_0}$. We will
show that $%
w_{1},w_{2}$ are in different Gleason parts of $A$.  In the 
definition in the introduction of $\M_{A_0}$,
take $\varepsilon=\frac{1}{m}$ for a natural number $m$, and choose the
corresponding $H_{m}=[h_{1}^{m}\cdots h_{n_{m}}^{m}],K_{m}=[k_{1}^{m}\cdots
k_{n_{m}}^{m}]\in A^{n_{m}}$ with $\left\langle H|K^{\ast}\right\rangle
=\sum_{i=1}^{n_{m}}h_{i}^{m}k_{i}^{m}=1$. 
By that definition we may assume
that $\Vert K_{m}(w)\Vert_{2}\leq c_{m}f(w)$ , and $\Vert H_{m}(w)\Vert
_{2}f(w)\leq c_{m}$, for all $w\in\Omega$. Here $c_{m}$ is a sequence of
real numbers decreasing to 1. Multiplying the functions\ $%
h_{i}^{m},k_{i}^{m} $ by constants with absolute value one we may also
assume that
\begin{equation}
h_{i}^{m}\left( w_{2}\right) \geq 0 \; \; 
\text{, for all }i=1,...,n_{m}.  \label{pos}
\end{equation}
Write $\Theta_{m}(w)=\frac{K_{m}(w)}{G\left( w\right) }$, and $\Pi
_{m}(w)=H_{m}(w)G(w)$. So $\Vert\Theta_{m}(w)\Vert_{2}\leq c_{m}$ and $%
\Vert\Pi_{m}(w)\Vert_{2}\leq c_{m}$. By the Cauchy-Schwarz Inequality we get
\begin{equation*}
1=|\sum_{i=1}^{n_{m}}h_{i}^{m}k_{i}^{m}|\leq\Vert\Theta_{m}(w)\Vert_{2}\Vert%
\Pi_{m}(w)\Vert_{2}\leq c_{m}\Vert\Theta_{m}(w)\Vert_{2}
\end{equation*}
Hence
\begin{equation*}
\frac{1}{c_{m}}\leq\Vert\Theta_{m}(w)\Vert_{2}\leq c_{m}\text{, for all }%
w\in\Omega.
\end{equation*}
Obviously, the same formula holds with $\Theta_{m}$ replaced by $\Pi_{m}$.
Expanding out the following square as an inner product and using (\ref{pos})
we get
\begin{equation*}
\Vert\Theta_{m}\left( w_{2}\right) -\Pi_{m}\left( w_{2}\right) \Vert
_{2}^{2}=\Vert\Theta_{m}\left( w_{2}\right) -\Pi_{m}\left( w_{2}\right)
^{\ast}\Vert_{2}^{2}\leq2c_{m}^{2}-2\rightarrow0.
\end{equation*}
Hence, by the Pythagorean Identity
\begin{equation*}
\Vert\Theta_{m}(w_{2})+\Pi_{m}(w_{2})\Vert_{2}\rightarrow2\quad\text{ as }%
m\rightarrow\infty.
\end{equation*}
Put $W_{n}=\frac{1}{2}\left( \Theta_{m}+\Pi_{m}\right) $. Since $%
h_{i}^{m},k_{i}^{m}\in A_{0}$ and $G\left( w_{1}\right) =1=-G\left(
w_{2}\right) $ we get
\begin{equation*}
\Vert W_{m}(w_{2})- W_{m}(w_{1})\Vert
_{2}=\Vert2\left( W_{m}\right) (w_{2})\Vert_{2}\rightarrow2\quad\text{ as }%
m\rightarrow\infty.
\end{equation*}

On the other hand, for any Euclidean vector $z$ of norm 1 we have 
$$ |(z.W_m)(w_{2}) - (z.W_m)(w_1)| \leq \left\| w_{2}-w_{1}\right\| 
\Vert z.W_m \Vert \leq \left\| w_{2}-w_{1}\right\| 
\sup \{ \Vert W_m(w) \Vert_2 : w \in \Omega \}\; ,$$
where $\left\| w_{2}-w_{1}\right\|$ is the norm of $w_{2}-w_{1}$ considered
as a functional on $A$.  
From the above it follows that $|(z.W_m)(w_{2}) - (z.W_m)(w_1)| 
 \leq \left\| w_{2}-w_{1}\right\| c_m$.  Thus
\begin{equation*}
\left\| W_{m}(w_{2})- W_{m}(w_{1})\right\|
_{2}\leq\left\| w_{2}-w_{1}\right\| c_m \; .
 \end{equation*}
Since $c_{m}\rightarrow1$, we get $\left\|
w_{2}-w_{1}\right\| =2$.  Thus
$w_{1},w_{2}$ lie in different Gleason parts of $%
A $.

To prove the other implication assume now that $w_{1},w_{2}$ are in
different Gleason parts of $A$. Choose a sequence of functions $r_{n}$
analytic in the disc $\bar{\Di}=\left\{ w\in \mathbb{C}:\left| w\right|
\leq 1\right\} $ with $r_{n}(1)=1,r_{n}(-1)=-1,$ and $1-1/n<|r_{n}(w)|<1+1/n$
for all $w\in \mathbb{D}$. Such functions can be found by taking a conformal
equivalence of $\mathbb{D}$ with a `smile shaped region' inside the annulus $%
1-1/2n<|w|<1+1/2n$, with the two tips of the smile at $-1$ and $1$. Since $%
w_{1},w_{2}$ are in different Gleason parts of $A$, for any $n\in \mathbb{N}$
large enough, there is an 
$a_{n}$ in $A$ such that $\left\| a_{n}\right\| \leq 1-\frac{1}{2n}$%
, $\left| (r_{n}\circ a_{n})\left( w_{1}\right) -\left( 1-\frac{1}{2n}\right)
\right| <\frac{1}{2n}$, $\left| (r_{n}\circ a_{n})\left( w_{2}\right) +\left(
1-\frac{1}{2n}\right) \right| <\frac{1}{2n}$. Also the norm, on $\left\{
a\in A:a\left( w_{1}\right) = 0 \right\} ,$ of the ``evaluation at $%
w_{2}$'' functional is equal to one.  Thus
there is a $\widetilde{a}\in A$ such that $\left\|
\widetilde{a}\right\| \leq \frac{1}{n},$ $\widetilde{a}\left( w_{1}\right)
=0$ and $\widetilde{a}\left( w_{2}\right) =
(r_{n}\circ a_{n})\left( w_{1}\right) +
(r_{n}\circ a_{n})\left(w_{2}\right)$. Put $b_{n}=r_{n}\circ a_{n}-%
\widetilde{a}\in A$. Notice that $b_{n}\left( w_{1}\right) G\left(
w_{1}\right) =b_{n}\left( w_{2}\right) G\left( w_{2}\right) $ so $b_{n}G\in
A_{0}$.  Moreover the spectrum of $b_{n}$ is contained inside the annulus $%
1-2/n<|w|<1+2/n$.  Hence $b_{n}$ is an invertible element of $A$, and
consequently $b_{n}G$ is an invertible elements of $A_{0}$. Thus $f=\lim
\left| b_{n}G\right| \in \bar{Q}_{A_{0}}$.  By Proposition \ref{ma} this
proves that $A_{0}f$ is a strong Morita equivalence bimodule.
\end{proof}

\vspace{4 mm}

\noindent {\bf Example.}  Let $A$ be a function algebra.
First take $\Omega = M_A$ and suppose that $A_0, G, f$ and 
$w_1,w_2$ are as in the last theorem.  Now, for any $\Omega$ on which
$A$ sits as a function algebra (such as $\partial A$), 
suppose further that there exists a function $h \in A$, such that 
$|h| = 1$ on $\Omega$ and 
$h(w_1) = - h(w_2)$.  We are not assuming $w_1,w_2 \in \Omega$ here, they 
are points in $M_A$.  Then $f = |G| = |hG| \in P(\Omega)$,
and $f^{-1} = |G^{-1}| = |\frac{h}{G}| \in P(\Omega)$.
Thus the submodule $A_0f$ of $C(\Omega)$ 
is a strong subequivalence $A_0-A_0$-bimodule, indeed a
unitary subequivalence bimodule.   However if 
$w_1, w_2$ are in the same Gleason part of $M_A$, then $A_0f$ 
is not a strong Morita equivalence $A_0-A_0$-bimodule.

For a very concrete example, let $A(\Di)$ for the disk algebra,
and set $\Omega = \T$, the unit circle.  Choose two points
$\alpha, \beta \in \Di$.  Then there exists an inner function $h$ in
$A(\Di)$ such that $h(\alpha) = -h(\beta)$.
For example, if $\alpha = -\beta = \frac{1}{2}$, let
$h(z) = z$.  Choose $G \in A(\Di)^{-1}$ such that
$G(\alpha) = - G(\beta)$.  Then the submodule $A_0|G|$ of 
$C(\T)$ is a strong 
subequivalence bimodule, which is not a strong Morita
equivalence $A_0-A_0$-bimodule.

\vspace{4 mm}

 Question: Suppose $Af$ is a
subequivalence bimodule where $\Omega$ is the maximal ideal
space of $A$.  Then is $Af$ a strong Morita equivalence bimodule?
In other words,
is $\G(M_A) \cap \G(M_A)^{-1} = \M^2$?   Recall that $P(M_A) \cap 
P(M_A)^{-1} = Q(M_A)$.  By the way, it follows easily from 
this latter fact, and the example above, that $P \cap P^{-1}$
is not a harmonic class.

\vspace{4 mm}

In the spirit of the calculation in the proof above, we end this
section by giving some alternative descriptions of 
$\mathcal{M}_{A}(\Omega )$.

\begin{proposition}
For a uniform algebra $A$ on a compact set $\Omega $, the 
class $\mathcal{M}_{A}(\Omega )$ coincides with
 the set of those functions $f$ $\in C(\Omega )^{+}$ for
which there exist a sequence of positive integers $\{n_{m}\}$ and a sequence
of $A$-tuples $H_{m},K_{m}\in A^{n_{m}}$, such that the following three
limits
\begin{itemize}
\item [(i)] $\Vert K_{m}(w)\Vert _{2}\rightarrow f(w)$ ,

\item [(ii)] $\Vert H_{m}(w)\Vert _{2}\rightarrow f(w)^{-1}$ , and

\item [(iii)]
 $\Vert K_{m}(w)-H_{m}(w)^{\ast }f(w)^{2}\Vert _{2}\rightarrow 0$
\end{itemize}
\noindent are valid uniformly over $w\in \Omega $. Here the `$\ast $'
represents the complex conjugate of the vector.
\end{proposition}

\begin{proof}
Suppose that $f \in \mathcal{M}_{A}(\Omega )$.
As in the proof of the previous theorem,
take $\epsilon =\frac{1}{m}$ for a natural number $m$, and
choose $K_{m}, H_{m}, c_m$ as in that proof.
Write $\Theta _{m}(w)=\frac{K_{m}(w)}{f(w)}$,
and $\Pi _{m}(w)=H_{m}(w)f(w)$. So $\Vert \Theta _{m}(w)\Vert _{2}\leq c_{m}$
and $\Vert \Pi _{m}(w)\Vert _{2}\leq c_{m}$. 
Again by Cauchy-Schwarz we get
\begin{equation*}
\frac{1}{c_{m}}\leq \Vert \Theta _{m}(w)\Vert _{2}\leq c_{m}
\end{equation*}
for all $w\in \Omega $; and the same formula holds for $\Pi _{m}$. 
By expanding out the following square as an inner product, we see
again that
\begin{equation*}
\Vert \Theta _{m}(w)-\Pi _{m}(w)^{\ast }\Vert ^{2}\leq
2c_{m}^{2}-2\rightarrow 0
\end{equation*}
uniformly as $m\rightarrow \infty $.

This gives the one direction of the proposition. 
However, the argument is
reversible until at the end we obtain that $H_{m}(w)\cdot
K_{m}(w)\rightarrow 1$ uniformly in $w\in \Omega $. Thus $%
b_{m}=H_{m}(w)\cdot K_{m}(w)$ is invertible in $A$, and $%
(b_{m}^{-1}H_{m}(w))\cdot K_{m}(w)=1$.
\end{proof}

\vspace{4 mm}

Notice that (i) and (ii) in the proposition
say exactly that $Af$ is a strong
subequivalence bimodule, or that $f^2 \in \G \cap \G^{-1}$.
Secondly, notice in fact that (iii) together
with either of (i) or (ii), implies the other condition.
However it is convenient to state it as such.

The last proposition may be
loosely phrased as saying that a function $f \in C(\Omega)^+$
is in $\M$ if and only if
there exists a sequence of $A$-tuples $K_m$,
such that (i) holds uniformly,
and $\frac{K_m(\cdot)^*}{\Vert K_m(\cdot) \Vert^2_2}$)
is `eventually, uniformly, an $A$-tuple'.

In view of (i) and (ii)
we may replace (iii) by the more appealing looking condition:
\begin{itemize}
\item [(iii)']  $\Vert K'_m(w) - H'_m(w)^*
 \Vert_2 \rightarrow 0$
\end{itemize}
uniformly on $\Omega$, where
$v' = \frac{v}{\Vert v \Vert}$ is the unit vector in the
same direction as $v$, for any vector $v$ in complex Euclidean
space.

\section{The Picard group.}
Some material in this
 section, and indeed in the rest of the 
paper,  requires  some technical knowledge
 of several papers.  The reader who desires
further background is directed to
\cite{Bnat} for a leisurely
introduction to our work, and
to \cite{BMP, BMN,Bhmo}
for more specific details.
The general reader is advised to simply read the main results
below.

We begin by discussing strong Morita equivalence of function 
algebras.  As we said earlier, we shall not give the general 
definition, but simply say that it involves
 a pair of bimodules $X$ and $Y$, called
\emph{equivalence bimodules}.
A well know property of (strong) Morita equivalence is that if $%
A$ and $B$ are unital, commutative and (strongly) Morita equivalent then $A$
is (isometrically) isomorphic to $B$. So we shall assume that $A=B$. However
this is not the end of the story, for the question remains as to which 
$A-A$-bimodules implement 
such a `self-equivalence' of $A$. The collection of
such bimodules, with two such bimodules identified if they are (completely)
isometrically $A-A$-isomorphic, is a group, with multiplication being the
Haagerup tensor product $\otimes _{hA}$ (see \cite{BMP}). 
We call this the \emph{strong
Picard group}, and write $Pic_{s}(A)$. For background on the Picard group,
see  any algebra text covering Morita equivalence (for example
\cite{F} Chapter 12), or \cite{BGR,DG} for a discussion of the C$^{\ast }-$%
algebraic version. For example we will show
below that for $A=A(%
\mathbb{D})$, the disk algebra, this Picard group is the direct product of
the Mobius group and the abelian group $C(\mathbb{T})/Re\;A$. 

In fact
for the  most part  we
will only consider a certain subgroup of $Pic_{s}(A)$, namely the
singly generated bimodules. We will
show that for a function algebra $A$ on a compact (Hausdorff) space $%
\Omega $, every such bimodule is essentially
of the form $Af$, where $f$ is a strictly
positive continuous function on $\Omega $.
Indeed $f$ may be chosen in $\M(\Omega)$, and by
\ref{gen} it follows that $f$ is unique up to the
coset $Q$.
 
Indeed, even in the non-singly generated case, if the second Cech
cohomology group of $\Omega $ vanishes, we shall see that every 
strong Morita equivalence $A-A$-bimodule is a finitely generated submodule of $C(\Omega )$.
Of course, the non-singly generated case is probably much more interesting,
but will no doubt require a much deeper analysis. Our intention here is
mainly to point out very clearly the features of the singly generated case.

As we said earlier, the multiplication on the Picard group of
$A$, is given by the module Haagerup tensor product $\otimes_{hA}$.
See \cite{BMP} for details about this tensor product.
However in the case that the modules are singly generated
(as just discussed), then this tensor product becomes rather trivial:

\begin{lemma}
\label{prod}
Suppose that $f \in \M(\Omega)$, and $g \in C(\Omega)^+$.
Then $Af \otimes_{hA} Ag
\cong A(f g)$ (completely) $A-A$-isometrically.
\end{lemma}

We omit the proof of this, which is a special case
of the more general Lemma \ref{sg}, which we prove later.

A standard type of
strong Morita equivalence bimodule for $A$  comes from
taking an isometric automorphism
$\theta : A \rightarrow A$, and defining the module
$A_\theta = A$, with the usual left module action,
and with right module action $b \cdot a = b \theta(a)$.
It is easy to check this is a strong Morita
equivalence bimodule for $A$.  The collection of
(equivalence classes of) this type of equivalence bimodule is
a subgroup of $Pic_s(A)$, which is isomorphic to
the group $Aut(A)$ of isometric automorphisms of $A$.  This,
in the case that $A$ is a function algebra,
corresponds to a group of homeomorphisms of the
maximal ideal space of $A$, which restrict to homeomorphisms
of the Shilov boundary of $A$.

More generally, if $X$ is a strong Morita equivalence $A-A$-bimodule,
and if $\theta \in Aut(A)$, then $X_\theta$ is also a
a strong Morita equivalence $A-A$-bimodule.  Here 
$X_\theta$ is $X$ but with
right action changed to $x \cdot a = x \theta(a)$.  One way to see 
this is to note that $X_\theta \cong X \otimes_{hA} A_\theta$.
 
An $A-A$-bimodule $X$ will be called `symmetric' if $a x = x a$ for all
$a \in A , x \in X$.

\begin{proposition} 
\label{semidi} For a function algebra $A$,
$Pic_s(A)$ is a semidirect product of $Aut(A)$ and the
subgroup of $Pic_s(A)$ consisting of symmetric
equivalence bimodules.  Thus,
every $V \in Pic_s(A)$ equals $X_\theta$, for a symmetric
$X \in Pic_s(A)$, and for some $\theta \in Aut(A)$.
\end{proposition}

\begin{proof}
This follows essentially as
in pure algebra (\cite{F} Chapter 12.18).
Suppose that $X$ is any strong Morita equivalence $A-A$-bimodule.
Then by the basic Morita theory,
any right $A$-module map $T : X \rightarrow X$ is simply left
 multiplication by a fixed 
element of $A$.  Indeed, via this identification, we have 
$A \cong CB_A(X)$ isometrically and 
as algebras (\cite{BMP} 4.1 and 4.2).  For fixed $a \in A$, 
the operator $x \mapsto xa$ on $X$, is a right $A$-module map,
with completely bounded norm $= \Vert a \Vert$.
Therefore by the above fact, there is a unique
$a' \in A$ such that $a'x = x a$ for all $x \in X$.  The map
$a \mapsto a'$ is then seen to be an isometric unital
automorphism $\theta$ of $A$.  In this
way we have defined a surjective group homomorphism 
$Pic_s(A) \rightarrow Aut(A)$.  This homomorphism has a 
1-sided inverse $Aut(A) \rightarrow Pic_s(A)$, namely
$\theta \rightarrow A_\theta$.   In this way, we see the
`semidirect product' statement.   Note that $X = 
(X_{\theta^{-1}})_\theta$,
and $X_{\theta^{-1}}$  is symmetric.
\end{proof}     

\vspace{4 mm}

More generally, 
if $X$ is a strong subequivalence $A-A$-bimodule in the sense of
\cite{BOMD}, then a similar argument works.  Since we will not
use this result here we will not give full details, but the 
idea goes as follows:  Suppose that $X$
 corresponds (using the notation of 
\cite{BOMD} \S 5) to a subcontext $(A,A,X,Y)$ of
a C$^*$-Morita context
$(\C,\D,W,Z)$, where $\C = C(\Omega)$.   It follows that 
$\D \cong \C$ isometrically.  
Now use the argument of the proof above, combined 
with Theorem 5.6 in \cite{BOMD}.

Thus we may assume henceforth that $X$ is symmetric. 
Then $X$
 dilates to a strong Morita equivalence 
$\C-\C$-bimodule $W = \C \otimes_{hA} X$.  
From \cite{Bhmo} Theorem 6.8, we know that
$W$ contains $X$ completely isometrically.
  It is helpful to consider the inclusion
$$
\left[ \begin{array}{ccl}
A & X \\
Y & A 
\end{array} \right] \subset
\left[ \begin{array}{ccl}
\C &  W \\
Z  & \C
\end{array} \right]
$$
of linking algebras.  Note that since $W \cong
\C \otimes_A X$ we have $wa = aw$ for all
$w \in W, a \in A$.  Similarly for
$Z \cong Y \otimes_A \C$ we have $za = az$.
Since $Z = W^*$, we have
$w a^* = (a w^*)^* = (w^* a)^* = a^* w$.  Therefore
$x w = w x$ for all $w \in W, x  
\in \C$.  Thus $W$ is a symmetric element of 
$Pic(C(\Omega))$, and consequently (see the appendix 
of \cite{Rae}, or \cite{Ri2,BGR,DG})
$W$ may be characterized as the space of sections
of a complex line bundle over $\Omega$.  

\begin{proposition}  Suppose that $A$ is a function algebra
on $\Omega$, and that $X$ is a symmetric subequivalence 
$A-A-$bimodule (or strong Morita equivalence 
bimodule). 
If every
complex line bundle over $\Omega$ is trivial, or 
equivalently, if
the Cech cohomology group $H^2(\Omega,\Z) = 0$, then
$W \cong C(\Omega)$ completely
$A-A$-isometrically.  Thus $X$ is $A-A$-isometric to
a finitely generated $A-A$-submodule of 
$C(\Omega)$.      
\end{proposition}

The hypothesis $H^2(\Omega,\Z) = 0$ applies, for example,
if $\Omega = \bar{\Di}$ or $\T$, or more generally, any
multiply connected region in the plane.

\begin{theorem}
\label{nice}  Suppose that $A$ is a function algebra on
$\Omega$.  Every singly generated strong Morita 
equivalence (resp. subequivalence)  
$A-A$-bimodule is completely
$A-A$-isometrically isomorphic
to $(Af)_\theta$, for some $\theta \in Aut(A), 
f \in C(\Omega)^+$.  The converse 
of this statement is also true, if $A$ is logmodular 
or if it is a logMorita algebra
(resp. convexly approximating in modulus).
\end{theorem}

\begin{proof}
We use the notation and facts established in the second paragraph
below the proof of Proposition \ref{semidi}.
If $X$ is  singly generated over $A$, then
so is $W$ over $\C$, and again it follows that
$W \cong \C$  completely
$A-A$-isometrically.  Here is one way to see this
(which adapts to the noncommutative case):
if $x$ is the generator,
then $x^*x$ is a strictly positive element in $\C$
(see for example \cite{BMP} theorem 7.13).  Thus
$x^*x$ is invertible, so that $f = |x|$
is invertible in $\C$.  Thus
$W = \C x \cong \C f = \C$, and
$X = Ax \cong Af$.

The `converse' assertion follows from the earlier correspondences
between the classes $\bar{Q}^+, \M$ and $\G$, and strong Morita
equivalence and subequivalence bimodules.
\end{proof}

\vspace{3 mm}

It can be shown that for $X$ a singly generated strong Morita
equivalence $A-A$-bimodule, 
the Morita `inverse bimodule'
$\tilde{X} \cong CB_A(X,A)$ may be identified
with $\{ c \in \C : c f \in A \} = A f^{-1}$.
That is, the strong Morita context (in the sense
of \cite{BMP}) associated with $X$
may be identified with  $(A,A,Af,Af^{-1})$.  The
latter is a subcontext (in the sense of \cite{BOMD})
of $(\C,\C,\C,\C)$.  Also, $f \in \M_A$.
These facts may be probably be shown
directly from what we have done in the proof, and \cite{Bhmo} Theorem
5.10 say.  However another proof of these facts is given in
Corollary \ref{Motat} in the next section.

We turn now to `rank one' modules.  We illustrate this concept 
first at the level of pure algebra.
Suppose that $A$ is a commutative unital algebra, and that
$X$ is a (purely algebraic)
Morita equivalence $A-A$-bimodule, with
`inverse bimodule' $Y$ (see \cite{F}).  As above, we say that
$X$ is symmetric if
$ax = xa$ for all $x \in X, a \in A$.  We shall say that
$X$ is {\em algebraically rank 1}, if $1 = (x',y') = [y'',x'']$ for some
$x', x'' \in X, y',y'' \in Y$; here $(\cdot)$ and
$[\cdot]$ are the Morita pairings \cite{F}.  It is
easy to show that this implies that $(x,y) = [y,x]$
for all $x \in X, y \in Y$.   Simple algebra shows that
a symmetric bimodule $X$ is
an algebraically rank 1 Morita equivalence $A-A$-bimodule
if and only if $X \cong A$ as $A$-modules.

Now suppose that 
$X$ is a
symmetric strong Morita equivalence $A-A$-bimodule.  We
 say that $X$ is {\em rank one}, if  $X$ is algebraically rank 1,
and  for any $\epsilon > 0$, 
the $x', y'$ above may be chosen with norms within $\epsilon$ of 1.
It follows, that if we define $T(x) = (x,y')$ on $X$, and 
$S(a) = ax'$ on $A$, then $S = T^{-1}$.  Also the norms
 $\Vert T \Vert_{cb} = \Vert T \Vert$, and $\Vert S \Vert_{cb}
= \Vert S \Vert$, are close to 1.   Thus it follows that 
$X \cong A$ almost completely $A$-isometrically, and we see
that $X$ is a $MIN$ space (that is, its operator space structure is
that of a subspace of a commutative C$^*-$algebra).

  We can therefore add to our earlier 
Corollary \ref{Motat21}, the
following:

\begin{corollary}
\label{Motat2}  Let $X$ be a Banach
$A$-module.  The following are equivalent:
\begin{itemize}
\item [(i)]  $X \cong A$ almost $A$-isometrically,
\item [(ii)]  $X$, with the symmetric bimodule action,
is a rank one  strong Morita equivalence $A-A$-bimodule,
\item [(iii)]  There exists $f \in \bar{Q}^+$ such that $X \cong Af$
 $A$-isometrically.
\end{itemize}
\end{corollary}

\begin{proof}  We just saw that $(ii) \Rightarrow (i)$.
As we saw in 
\ref{Motat21}, condition (i) is equivalent to (iii).
Proposition \ref{ma} (c) shows that (iii) implies (ii).
\end{proof}

\vspace{2 mm}

Thus the new definition of rank one
is equivalent to our earlier
definition (\ref{prec} (i))
of $Af$ being a rank one strong Morita equivalence
bimodule.

\begin{theorem}
\label{gen}
\begin{itemize}
\item [(i)]  If $f_1, f_2 \in \M$
then $Af_1 \cong Af_2$ $A$-isometrically
if and only if $f_1 = h f_2$
for some $h \in Q$.
\item [(ii)]  The map $[f] \mapsto Af$ gives an injective
group homomorphism $\M' = \M/Q \rightarrow Pic_s(A)$.
The range of this homomorphism consists of
all the topologically
singly generated symmetric elements of $Pic_s(A)$.
Thus $Pic_s(A)$ contains as a subgroup,
the direct product of $\M'$ and $Aut(A)$.
\item [(iii)]   If one restricts the
 group homomorphism in (ii) to
$Q' = \bar{Q}^+/Q$, then its range consists
of the rank 1 strong Morita equivalence $A-A$-bimodules.

\item [(iv)]  If $f_1, f_2 \in \M$
then $Af_1 \cong Af_2$ almost $A$-isometrically
if and only if $f_1 = h f_2$,
for some $h \in \bar{Q}^+$.
\end{itemize}
\end{theorem}

\begin{proof}  By the `harmonicity' associated with the
$\M$-class, we can assume that $f_1, f_2 \in C(\partial A)^+$.
Then (i) follows immediately from Corollary \ref{wasla}.
However here is another
argument, which adapts immediately to give (iv) too.
  Suppose that $Af_1 \cong Af_2$ $A$-isometrically,
where $f_1, f_2 \in \M$.  Then
by \ref{prod} we have
$A \cong Af_1^{-1} \otimes_{hA} Af_1 \cong
Af_1^{-1} \otimes_{hA} Af_2 \cong A(f_2f_1^{-1})$.
This implies that $f_2f_1^{-1} \in Q$ by the last theorem.
Thus $[f_1] = [f_2]$.  This also gives (ii), in view of
\ref{prod}.
A similar argument to (i) also proves (iv).
\end{proof}

\vspace{5 mm}

For many common function algebras, every strong Morita
equivalence $A-A$-bimodule $X$ is singly generated.
We shall illustrate this for $A = A(\Di)$.  Such an $X$
is algebraically finitely generated and projective 
(\cite{F} 12.7).
Then by
Theorem 3.6 in \cite{Bhmo}, $X$ is completely 
boundedly $A$-isomorphic to a closed
$A$-complemented submodule of $A^{(n)}$ .  The associated
projection $P : R_n(A) \rightarrow R_n(A)$ may be thought
of as an analytic projection valued function 
$f_P : \Di \rightarrow
M_n$, and since the disk is contractible, 
$f_P$ is homotopic to the
constant function $f_P(0)$.  Thus (see \cite{Bla} 4.3.3)
$P$ is similar to a constant projection in $M_n$.  We 
originally heard this last argument from P. Muhly. 
Hence $X \cong
A^{(m)}$ algebraically, and if $m > 1$ then $W = 
\C \otimes_A X \cong \C^{(m)}$ algebraically,
 which is impossible.   So $X$ is singly generated.  
Putting this together with Theorem \ref{gen} (ii), and the fact 
that $A(\Di)$ is a Dirichlet algebra (so that 
$\bar{Q}^+ = \M = C(\T)^+$), we have: 

\begin{corollary}
\label{name}  The Picard group of $A = A(\Di)$, 
is the direct product of the Mobius group
and the abelian group $C_{{\mathbb R}}(\T)/ Re \; A $.
\end{corollary}

\noindent {\bf Remark.}  We have proved elsewhere
that the strong Picard group of
$A$ is isomorphic to the group of category equivalences
of $_AOMOD$ with itself, where $_AOMOD$ is the category of left
operator modules over $A$.

\vspace{3 mm}

\begin{center}
Part C.
\end{center}

\vspace{3 mm}

\section{Rigged modules over function 
algebras.}

A (left) $A$-Hilbertian module is a left operator module
$X$ over $A$, such that there exists a net of natural
numbers $n_\alpha$ , and completely contractive
$A$-module maps $\phi_\alpha
: X \rightarrow R_{n_\alpha}(A)$ and $\psi_\alpha :
R_{n_\alpha}(A) \rightarrow X$, 
such that $\psi_\alpha \circ \phi_\alpha
(x) \rightarrow x$, for all $x \in X$.   
Here $R_n(A)$ is $A^{(n)}$, viewed as an operator space
by considering it as the first row of $M_n(A)$.
The name 
`Hilbertian' is due to V. Paulsen. 
  An $A$-rigged 
module is an $A$-Hilbertian module with the additional
property that for all $\beta$ we have
 $\phi_\beta \psi_m \circ \phi_m
\rightarrow \phi_\beta$ in cb-norm.   
If in addition,
 $X$ is singly generated, it follows that we can take 
the net in the definition of 
$A$-Hilbertian, to be a sequence $n_m,  \phi_m, \psi_m$,
$m \in \N$.  We will write $e_m = \psi_m \circ \phi_m$,
a completely contractive module map $X \rightarrow X$.      
 We will also say that
$X$ is {\em rank 1 Hilbertian} if we can take
 $n_m = 1$ for all $m \in \N$.

Our main purpose in this section is to 
show that a singly generated operator
module is $A$-Hilbertian if and only if it is $A$-rigged,
and to attempt to thoroughly understand
such modules. 
It is clear from the definitions, that
every strong Morita equivalence  $A-A$-bimodule is
a left $A$-rigged module, and every left $A$-rigged module
is an $A$-Hilbertian module.  For C$^*-$algebras, we proved 
in \cite{Bna,BMP} that the converse is true.  

If $E$ is a closed subset of $\Omega$, we will write $J_E$
for the ideal $\{ f \in A : f(x) = 0 \; \text{for all} \;
x \in E \}$  of $A$.    The following result is important for us:  

\begin{theorem} \cite{Sm,HWW,ER}
\label{Smith} Let $A$ be a function algebra on compact
$\Omega$.  If $J$ is an ideal
in a function algebra $A$, then
the following are equivalent:
\begin{itemize}
\item [(i)]  $J$ has a contractive approximate identity
\item [(ii)]  $J$ has a bounded approximate identity
\item [(iii)]  $J$ is an M-ideal of $A$,
\item [(iv)]  $J = J_E$, for a p-set $E$ for $A$ in $\Omega$.  
\end{itemize}
\end{theorem}

We refer the reader to \cite{Gam} for details on p-sets
and peak sets.  We allow $\emptyset$ as a p-set.
A p-set is an intersection of peak sets.  We will not define the 
term `M-ideal' here, and it will not play a role.

In (iv), it is clear that the $E$ is unique.  Also
$F = E \cap \partial A$ is a p-set for $A$ on the
Shilov boundary, and $J_F = J_E$.
 
It is easy then to see from the definition of a rigged module
above that we have:

\begin{corollary}
\label{co}
If $E$ is a p-set for a function algebra $A$,
then the  ideal $J_E$, is an $A$-rigged module.
\end{corollary}

\begin{proof}
Let $e_\alpha$ be a c.a.i. for $J_E$, and define
$\psi_\alpha = \phi_\alpha$ to be multiplication by
$e_\alpha$.  These
satisfy the requirements for a rigged module.
\end{proof}

\vspace{5 mm}

The proof of the following corollary
requires some technical knowledge of
rigged and C$^*-$modules.

\begin{lemma}
\label{rig}  Suppose that 
$A \subset
C(\Omega)$ is a function algebra on a compact space $\Omega$,
and that $X$ is a 
singly generated left $A$-Hilbertian module.
Then there is a nonnegative
continuous function $f$ on $\Omega$, 
such that $X \cong (Af)^{\bar{}}$ $A$-isometrically.   
\end{lemma}

\begin{proof}  As observed in
\cite{Bna} \S 7, $W = C(\Omega) \otimes_{hA} X$ is a 
(singly generated) left C$^*-$module over $C(\Omega)$.  
Also $X$ may be regarded as an $A$-submodule of $W$
(in the obvious way).
Let $f$ be the single generator of $X$ and $W$.
Suppose that $I$ is the ideal in $\C = C(\Omega)$
generated by the range of the inner product.  
Now $W$ is a full C$^*-$module over $I$, and so 
$IW$ is dense in $W$, since C$^*-$modules are automatically
nondegenerate.  Since $\C f$ is dense in $W$, $I \C = I$, 
and $IW$ is dense in $W$,
we see that $If$ is dense in $W$.  Hence $W$ is singly generated
by $f$ over $I$.  The obvious map $F: I \rightarrow If \subset W$
is adjointable and $F, F^*$ have dense range, so
by the basic theory of C$^*-$modules (see e.g.
\cite{La} Prop. 3.8), $W \cong I$, $I$-isometrically.  Hence
$W \cong I$ , $\C$-isometrically.  
Since $W$ is singly generated by $f$, 
$|f|$ is a strictly positive
element in $I$,
by the argument before Theorem \ref{nice}
above.  Of course $|f|$ is
not strictly positive on $\Omega$, unless $I = A$.    
Clearly $X$ is $A$-isometric to the closure of the submodule
$A|f|$  of $I$.    
\end{proof}

It will be useful to have the following

\begin{lemma} \label{mmag} Suppose that $f \in C(\Omega)_+$, and  
that $K \in A^{(n)}$.  Then $\Vert a(w) K(w) \Vert_2
\leq \Vert af \Vert_{\Omega}$ for all $a \in A$ and 
$w \in \Omega$, if and only if
$\Vert K(w) \Vert_2 \leq f(w)$ for all $w \in \partial A$.
\end{lemma}

\begin{proof}  $(\Leftarrow)$:   For 
$w \in \partial A$,  we have 
$$\Vert a(w)K(w) \Vert_2 \leq |a(w)f(w)| \leq \Vert af 
\Vert_{\Omega} \; \;.$$
Since $aK$ achieves its maximum modulus on $\partial A$,
we have proved this direction.

$(\Rightarrow)$:  If $n=1$, then this is well known, 
following by the usual Choquet boundary point
argument (such as we've seen, for example, in 
\ref{NEED}).  For general $n$, fix
 $z \in \Co^n_1$.  Then we have
$|a(w)(z.K(w))| \leq \Vert af \Vert_{\Omega}$, for all 
$w \in \Omega$.   Thus by the $n=1$ case, 
$|z.K(w)| \leq f(w)$ for all $w \in \partial A$,
which gives what we need.
 \end{proof}

\begin{corollary}
\label{Hisr}
For a singly generated operator
module $X$ over a function algebra
$A$, $X$ is $A$-Hilbertian  if  and only if $X$ is $A$-rigged.
\end{corollary}

\begin{proof}
Suppose that $X$ is $A$-Hilbertian.
We may suppose, by \ref{rig},
that $X = (Af)^{\bar{}}$ in $C(\Omega)$.
We will take $\Omega$ to be the Shilov boundary of $A$.
We use the notation of the beginning of the section,
but assume, as we may, that the norms of $\phi_m ,
\psi_m$ are strictly less than 1.
The map $\psi_n$ may be written as
$[a_1, \cdots , a_k] \mapsto \sum_i a_i x_i$ for some
$x_i \in X$,
and without loss of generality, we can assume that
$x_i = h_i f$, with $h_i \in A$.  Thus $\psi_n$ may be associated 
with an $A$-tuple $H_m$, and without loss of generality,
$f(w) \Vert H_m(w) \Vert_2 \leq 1$ for all $w \in \Omega$.
See \cite{Bhmo} Prop. 2.5(ii).  

Also, $\phi_m$ is 
completely determined by its
action on $f$, so we can associate $\phi_m$ with a unique
$A$-tuple $K_m$. 
Now it is easily seen that  
without loss of generality,
$e_m$ may be regarded as (multiplying by) 
the element $H_m(w) . K_m(w)$ of $A$.  We have
$f e_m \rightarrow f$ uniformly.   

Next we note that by Lemma \ref{mmag}, the identity 
$\Vert \phi_n(af) \Vert \leq 
\Vert af \Vert_{\Omega}$ for all 
$a \in A$,  implies that $\Vert K_m(w) \Vert_2
\leq f(w)$ for all $w \in \Omega$.  
Thus it follows that for $a \in A$, we have
$$ \Vert e_m a \Vert_\Omega \leq C_m \Vert af \Vert_\Omega$$
for a constant $C_m$ that does not depend on $a$.

Set $\phi'_m = \phi_m e_m, \psi'_m = \psi_m$.
These give new factorization nets, but now we have that 
$$\Vert \phi'_m - \phi'_m 
\psi'_n \phi'_n \Vert_{cb} = 
\Vert \phi_m e_m (1 -  e_n^2) \Vert_{cb} 
\leq \Vert e_m(1 - e_n^2) \Vert_\Omega 
\leq C_m \Vert f(1 - e_n^2) \Vert_\Omega
\rightarrow 0$$
 as $n \rightarrow \infty$.
This says that $X$ is a rigged module.
The converse is trivial.   
\end{proof}

\vspace{5 mm}

\begin{definition}
\label{arr}
If $E$ is a closed subset of $\Omega$, we define 
$\R_E(\Omega)$ to be the set of functions $f \in C(\Omega)$ 
which vanish exactly on the set $E$, and 
for which there exists two sequences $H_n$ and $K_n$
 of $A$-tuples, such that 
\begin{itemize}
\item [(i)] $ \Vert K_n(w) \Vert_2 \leq f(w)
\leq \Vert H_n(w) \Vert_2^{-1}$ for all
$w \in \Omega \; $ (interpreting $\frac{1}{0} = \infty)$,
\item [(ii)]  $e_m(w) = H_m(w) . K_m(w) \rightarrow 1$ uniformly
on compact subsets $C$ of $\Omega \setminus E$ . \\

\noindent We will write $\R(\Omega)$ for 
the combined collection of all the 
$\R_E(\Omega)$ classes.
\end{itemize}
\end{definition}

If we write $g = \log f$, then $f \in \R_E(\Omega)$ if and
only if $g$ is finite precisely on $E$, and there exist
$H_n, K_n$ as above, satisfying (ii) and
$$(i') \; \; \log \Vert  K_n(w) \Vert_2 \leq g(w) \leq
- \log \Vert H_n(w) \Vert_2 \; \; . $$
It follows from these, that $\log \Vert  K_n(w) \Vert_2$ 
and $- \log \Vert H_n(w) \Vert_2$ converge uniformly to
$g$, on compact subsets of $\Omega \setminus E$ (see Lemma 
\ref{sil} below).
Thus we see that $g$ is an upper and lower envelope
of `sub- and superharmonic' functions.  See Proposition
 \ref{sil2}
for more on this.  This definition is therefore somewhat
reminiscent of the Perron process of solving the Dirichlet 
problem.

\begin{proposition}  
\label{req} $\R_{\emptyset}(\Omega)  = \M(\Omega)$.
\end{proposition}

We leave this as an exercise for the reader.
 
\begin{corollary}
\label{H}  Let $A$ be a function algebra on $\Omega$.  Then
\begin{itemize}
\item [(i)]  If $\R_E(\Omega) \neq \emptyset$ then $E$ is a 
peak set for
$A$ in $\Omega$.   
\item [(ii)]  If $f \in \R_E(\Omega)$, 
$(Af)^{\bar{}}$ is a rigged module. 
\item [(ii)]  If $\Omega = \partial A$ and  $f \in C(\Omega)_+$,
then $(Af)^{\bar{}}$ is a
rigged module if and only if 
$f \in \R_E(\Omega)$ for some peak set $E$.
\end{itemize}
\end{corollary}

\begin{proof}  If one studies the proof of Corollary \ref{Hisr},
one sees that the ideas there yield that if 
$f \in \R_E(\Omega)$ then
$(Af)^{\bar{}}$ is a rigged module.  If $E = f^{-1}(0)$, then
 $K_n$ vanishes on $E$.  By (ii) of Definition
\ref{arr}, the functions
$e_m(w) = H_m(w) . K_m(w)$ form a c.a.i. for 
$C_0(\Omega \setminus E)$, and also for $J_E$.
We deduce from Theorem
 \ref{Smith} that $E$ is a p-set.  Since $f$ is strictly positive
on $\Omega \setminus  E$, it follows that $E$ is a $G_\delta$,
from which we deduce, using \cite{Gam} II.12.1, that $E
$ is a peak set.  This gives (i) and (ii).

If $(Af)^{\bar{}}$ is a 
rigged module then the ideas of the proof of Corollary \ref{Hisr}
show that
(i) of Definition \ref{arr} holds 
for $w \in \partial A$. 
Set $E = f^{-1}(0)$, then $K_m$ vanishes on the set
$F = E \cap \partial A$.  We saw in that Corollary that 
$f e_m \rightarrow f$ uniformly on $\partial A$. 
By the Stone-Weierstrass theorem $C_0(\partial A \setminus F) f$
is dense in $C_0(\partial A \setminus F)$.  Thus the functions
$e_m(w) = H_m(w) . K_m(w)$, which are in $J_F$,
form a c.a.i. for $C_0(\partial A \setminus F)$.  
Thus by Urysohn's lemma,
the restriction of $f$ to $\partial A$ is in $\R_F(\Omega)$.
\end{proof}

\begin{corollary} \label{ops}  The singly generated rigged
modules over a function algebra $A$, are exactly (up to 
$A$-isometric isomorphism), the modules of form $(Af)^{\bar{}}$
for $f \in \R_E(\Omega)$.
\end{corollary}
 
If a subset $E \subset \Omega$ is a peak set for $A$ in $\Omega$,
then clearly  $F = E \cap \partial A$ is a peak set for $A$ in
$\partial A$ .  Conversely any peak set $F$ for $A$ in
$\partial A$  may be written $F = E \cap \partial A$ for a unique
peak set $E$ for $A$ in $\Omega$.  We will therefore sometimes
be sloppy, and write $\R_F(\Omega)$ or $\R_D(\Omega)$ for
$\R_E(\Omega)$ , where $D$ is the unique peak set for $A$ in $M_A$
with $D \cap \Omega = E$.

\begin{lemma}
\label{sil}
If $f, H_n , K_n$  is as in the definition of $\R_E(\Omega)$ , 
then $\Vert K_n(w) \Vert_2 \rightarrow 
f(w)$ uniformly on $\Omega$, and $\Vert H_n(w) \Vert_2
\rightarrow f(w)^{-1}$ uniformly on compact subsets of
$\Omega \setminus E$.  Also, $\Vert a f \Vert_{\partial A}
= \Vert a f \Vert_\Omega$ for all $a \in A$.  
\end{lemma}

\begin{proof}  If $C$ is a compact subset of $\Omega \setminus E$,
and $\epsilon > 0$ is given, then 
$$(1 - \epsilon) f(w) 
\leq f(w) |H_m(w) . K_m(w)| \leq \Vert K_m(w) \Vert_2 \leq f(w)$$
 uniformly for 
$w \in C$ and $m$ sufficiently large.  From this one easily sees
that $\Vert K_n(w) \Vert_2 \rightarrow
f(w)$ uniformly on $\Omega$, since $f$ is bounded away from zero 
on $C$.  The second statement is similar.  Finally,
for $w \in \Omega$,  we see that
$$|a(w)f(w)| \leq \sup_m \Vert a(w)K_m(w) \Vert_2
\leq \sup_m \sup \{ \Vert a(x) K_m(x) \Vert_2 : x \in
\partial A \}
\leq \Vert a f \Vert_{\partial A} \; \; .$$
\end{proof}

If $H$ is an $A$-tuple, then it will be useful to think
of $\log \Vert H(w) \Vert_2$ as a {\em subharmonic} function.
The logarithm of
a function $f \in \R_E$ on the other hand should be thought
of as being {\em harmonic}, as we mentioned briefly before.   
The following few results begin to justify these assertions:

\begin{proposition}
\label{sil2}
Let $A$ be a uniform algebra on compact $\Omega$, and suppose that 
$f \in  \R_E(\Omega)$.  Then:
\begin{itemize}
\item [(i)]   $f_{|_{\partial A}} \in \R_{E \cap \partial A}(\partial A)$.
\item [(ii)]  If  $g = \log f$, then $g$ achieves
its maximum and minimum on $\partial A$.
\item [(iii)]  If there is a domain $R \subset \Co^n$, and 
an inclusion $R \subset \Omega$, such that all functions in $A$ 
are analytic functions on $R$, then $g = \log f$ is
harmonic (in the usual sense on $R$) whenever it is finite
(that is, on $R \setminus E$).
\end{itemize}   
\end{proposition}

It is possible that in (iii) above, $R \setminus E = R$ 
automatically,
if $R$ is a (connected) domain as in (iii).  This is the case if 
$A$ is the disk algebra (see comments after Example \ref{ex4}).
We have not checked this in general though.

\begin{proof} (i) is obvious.  (ii):  Let $K_n, H_n$ be $A$-tuples 
as in Definition \ref{arr}.  For any unit vector $z$ in the complex
Euclidean space of the appropriate dimension, and any $w \in \Omega$,
 we have 
$$ |K_n(w).z| \leq  \Vert K_n(\cdot).z \Vert_{\partial A}
\leq \Vert f \Vert_{\partial A} \; \; . $$ 
Thus $\log \Vert K_n(w) \Vert_2 \leq \log \Vert f \Vert_{\partial A}$. 
Letting $n \rightarrow \infty$,  gives $g(w) \leq  \sup_{\partial A} g$.
To get the other inequality, we may assume that $g$ is bounded 
below.  Thus $E = \emptyset$.
As above, we obtain $\Vert H_n(w) \Vert_2 \leq
\sup_{\partial A} |f|^{-1}$.   Hence $- \log \Vert H_n(w) \Vert_2
\geq \inf_{\partial A} g$.  Now let $n \rightarrow \infty$.    

(iii):  
A function is harmonic on a domain in $\Co^n$ if 
and only if it satisfies the Mean Value Principle.   Let
$H_n, K_n, z$ be as in (ii), and fix $w_0 \in R \setminus E$.
Since $K_n(w).z$ is analytic for $w \in R$, we have that 
$\log |K_n(w).z|$ is subharmonic on $R$.
Thus, for any ball $B$ center 
$w_0$ in $R$, we have
$$\log |K_n(w_0).z| \leq \frac{1}{m(B)} \int_B \log 
|K_n(w).z| \leq \frac{1}{m(B)} \int_B \log
\Vert K_n(w) \Vert_2 \leq \frac{1}{m(B)} \int_B \log f(w) \; \; .
$$
Since this is true for all such $z$ we have
$\log \Vert K_n(w_0) \Vert_2 \leq \frac{1}{m(B)} \int_B \log f(w)$.
Taking the limit as $n \rightarrow \infty$ gives
$\log f(w_0) \leq \frac{1}{m(B)} \int_B \log f(w)$.  A similar 
argument using $H_n$ gives the other direction of the  Mean Value
Property.
\end{proof}

We state the following obvious fact since it will be referred to
several times:  

\begin{lemma}
\label{sil3}  Let $A$ be a uniform algebra on compact $\Omega$, and suppose that
$H$ and $K$ are $A$-tuples, such that
$\Vert H(w) \Vert_2 \Vert K(w) \Vert_2 \leq 1$
for all $w \in \Omega$, or all $w \in \partial A$.
Then the same inequality holds for
all $w \in M_A$.
\end{lemma}

\begin{proof}  Let $z,y$ be vectors in complex Euclidean space.
Then we have $|(z \cdot H(w))(y \cdot K(w))| \leq
1$ for all $w \in \Omega$, and consequently
for all $w \in M_A$.
\end{proof}

\begin{theorem}
Suppose 
that $E$ is a peak set for $A$ , and that 
$f \in \R_{E}(\Omega)$.  Let $D$ be the unique peak set
in $M_A$ with $D \cap \Omega = E$.  Then there 
is a unique function $\tilde{f}
\in \R_{D}(M_A)$ such that
$\tilde{f}$ restricted to $\Omega$ equals $f$.
Also, $\Vert af \Vert_\Omega = \Vert a \tilde{f} \Vert_{M_A}$
for all $a \in A$, so that $(Af)^{\bar{}}$ in 
$C(\Omega)$ is $A$-isometric to
$(A\tilde{f} )^{\bar{}}$ in $C(M_A)$.
\end{theorem}

\begin{proof}  
As we saw above,
condition (ii) of definition \ref{arr} 
is equivalent to saying that the functions
$e_m(w) = H_m(w) . K_m(w)$, which are in $J_E$, form a 
c.a.i. for $C_0(\Omega \setminus E)$.  However,
$J_E = J_D $, so that $e_m$ is a c.a.i. for $J_D$.   
Hence by \cite{Gam} II.12.5 and II.12.7, taking the 
function $p$ there 
to be a continuous strictly positive function which is
$1$ on $D$ and $< \epsilon$ on a compact subset $C$
of $M_A \setminus D$, we see that 
$ \; e_m(w) \rightarrow 1$ uniformly on $C$.
 
Now by lemma \ref{sil3},
we have
$\vert e_m(w) \vert \leq \Vert H_m(w) \Vert_2 \Vert K_m(w) \Vert_2
\leq 1 $
for all $w \in \Omega$.
This then implies that
$\lim_m \Vert H_m(w) \Vert_2 \Vert K_m(w) \Vert_2
= 1$ uniformly on any compact subset $C$
 of $M_A \setminus D$.  
Since $\Vert H_m(w) \Vert_2 $ and $\Vert K_m(w) \Vert$ are uniformly 
bounded above on $C$, 
we deduce that they are also uniformly bounded 
away from zero on $C$.   By Lemma \ref{sil3}, we have 
\begin{equation} 
 \Vert H_m(w) \Vert_2 \Vert K_n(w) \Vert_2 \leq 1 \; \; 
 \text{for all } w \in M_A \; .  \label{star} 
\end{equation}
  Thus we have  $\Vert H_m(w) \Vert_2 
\Vert H_n(w) \Vert^{-1}_2 \leq 1+\epsilon$ uniformly on
$C$, for $m,n$ large enough.  Thus the sequence
$\Vert H_n(w)\Vert_2$, and by symmetry the sequence $\Vert K_n(w) 
\Vert_2$,
are uniformly Cauchy on $C$.  
Let $\tilde{f}$ be the 
uniform limit of $\Vert K_n(w) \Vert_2$ on $C$.  Varying
over all compact $C$ gives a well defined continuous $\tilde{f}$
on $M_A \setminus D$.  Clearly $\tilde{f}$ extends $f$.  

We next show that $\tilde{f}
\in C_0(M_A \setminus D)$.  Let $\epsilon > 0$ be 
given, and let $C = \{ w \in \partial A : f(w) \geq \epsilon \}$.  
For $\gamma > 0$ to be determined, choose (by \cite{Gam} II.12 again) 
$a \in A$ with $a \equiv 1$ on $D$, $\Vert a \Vert
\leq 1 + \gamma$, and $|a| < \gamma$ on $C$.  For 
$x \in \partial A, m \in \N$
and any Euclidean vector $z$ of norm
1, we have $|(z \cdot K_m(x))a(x)| \leq 2 \epsilon$ , if
$\gamma$ is smaller than a certain constant which depends only
on $\Vert f \Vert_{\partial A}$ and $\epsilon$.   Hence
$|(z \cdot K_m(w))a(w)| \leq 2 \epsilon$, for all $w \in M_A$. 
 Hence $\Vert K_m(w) \Vert_2 |a(w)| \leq 2 \epsilon$.
Letting $m \rightarrow \infty$ we see that 
$\tilde{f}(w) |a(w)| \leq 2 \epsilon$.  In particular, for 
$w \in U = |a|^{-1}((1-\epsilon,\infty))$, we have 
$\tilde{f}(w) \leq \frac{2 \epsilon}{1-\epsilon}$. 
Thus indeed $\tilde{f} \in C_0(M_A \setminus D)$. 

Define $\tilde{f}$ to be zero on $D$.  
Clearly 
$\Vert K_m(w) \Vert_2 \rightarrow \tilde{f}(w)$ for all
$w \in M_A$.  
For $w \notin D$ we obtain from (\ref{star}), that
$\Vert H_m(w) \Vert \tilde{f}(w) \leq 1$, and also  
$\Vert K_m(w) \Vert_2 \leq \tilde{f}(w)$.   
Thus $\tilde{f} \in \R_E(\M_A)$.  
 
Finally, for the uniqueness, we suppose that 
$f_1$ with $H^1_m, K^1_m$, and $f_2$ with $H^2_m, K^2_m$, 
both fulfill the definition of $\R_E(M_A)$.  If
 $f_1(x) = f_2(x)$ for all $x \in \partial A$, then
$\Vert H^1_m(w) \Vert_2 \Vert K^2_m(w) \Vert_2 \leq 1$
on $\partial A$, and hence, by Lemma \ref{sil3}, on $M_A $.
Hence $f_2(w) \leq f_1(w)$ for any $w \notin D$.  This obviously
implies what we want, by symmetry.
\end{proof}

\vspace{5 mm}

In view of the previous result we may simply write
$\R_E$ , for $\R_E(\Omega)$, if we wish.  
As remarked earlier, we may 
switch $E$ for the corresponding peak set in $M_A$ or 
$\partial A$.

\begin{example}
\label{ex4}
\end{example}   The nontrivial
p-sets $E$ for the disk algebra
$A(\Di)$ coincide with the peak sets, and they
 are exactly the closed subsets of
$\T$ of Lebesgue measure $0$ (see \cite{HWW} for example).
By the well known version of
Beurlings theorem for the disk algebra,
$J_E$ is a singly generated $A(\Di)$-module.
Hence, for $E \subset \T$ with $|E| = 0$, we have
 by Corollary \ref{co} that  $J_E$ is an example of a singly
generated rigged module over $A(\Di)$.  Indeed
every closed ideal $I$ of $A(\Di)$ is isometrically isomorphic
to some $J_E$, and consequently is an $A$-rigged module.
This is because, by Beurlings theorem, $I = J_E g$ for a
fixed inner function $g$ and peak set $E$.

Next we find all singly generated rigged modules over
$A(\Di)$.  They all turn out to be `rank one Hermitian'.

\begin{theorem}
\label{ad}
Suppose that $f \in C(\T)_+$.  The following are
equivalent:
\begin{itemize}
\item [(i)]  $\log f$ is integrable on $\T$,
\item [(ii)]  $f = |\phi|$ for
 a function $\phi \in H^\infty \setminus
\{0 \}$,
\item [(iii)]  $(A(\Di) f)^{\bar{}}$ is a rigged module over $A(\Di)$,
\item [(iv)]  $f \in \R(\T)$.
\end{itemize}
Moreover, every singly generated rigged module over $A(\Di)
$ is
$A(\Di)$-isometric to one of the form in (iii).
\end{theorem}

\begin{proof}  The equivalence of (i) and (ii) is classical
(\cite{Ho2} p. 53).

Suppose that $X$ is a singly generated rigged module over $
A(\Di)$.
By Corollary \ref{ops}
and the fact mentioned in the previous paragraph,
we have  $X \cong (A(\Di)f)^{\bar{}}$
$A$-isometrically, where
$E$ is a subset of $\T$ of Lebesgue measure $0$,
and $f \in \R_E(\T)$.
With the earlier notation, we have
 $\Vert K_m(w)  \Vert_2 \leq f(w)$, and
$\Vert K_m(w)  \Vert_2 \rightarrow f(w)$,
for all $w \in \T$.  This implies that there is a nonzero
function $K \in A(\Di)$, such that $|K| \leq f$ on
$\T$.  Then $\log |K| \leq \log f$, a.e.
on $\T$.  Since $\log |K|$ is integrable
 on $\T$, so is $\log f$ .

Conversely, let $\phi \in H^\infty$, with $f = |\phi|$ 
continuous on
$\bar{\Di}$.   Let $E$ be the subset of $\T$ on which $\phi$
vanishes, which is a a closed subset of measure $0$.
Let $w = \log |\phi|$, then $w$ is integrable.
We choose a function $k_1$ on
$\T$ such that $w - \epsilon \leq k_1 \leq w$, and
such that $k_1$ is continuously differentiable wherever
it is finite.  Indeed, one may assume that $k_1$ lies in
a thin strip about $w - \frac{\epsilon}{2}$.  We define
$$k(z) = \exp \left( \frac{1}{2\pi} \int_{-\pi}^{\pi}
\frac{e^{i \theta} + z}{e^{i \theta} - z} \; k_1(\theta) \;
 d \theta
\right) \; \; .
$$
Then $k \in A(\Di)$, and $\vert k \vert \leq \vert
 \phi \vert$
on $\T$.

Next, choose an open subset $U$ of $\T$ containing $E$, such
that $\int_U |k_1| < \epsilon$, and such that
$k_1 < -1$ on $U$.  We may suppose that
$U$ is a finite collection of disjoint open intervals, of total
combined length $< \epsilon$.
Choose $k_2 = k_1 + \epsilon$ outside $U$ .  On $U$ we define
$k_2$ so that $k_2$ lies between $0$ and $w$ on $U$, and so
 that
$k_2$ is finite and continuously differentiable
on all of
$\T$.  It is not hard to see that this is possible.  Then we have
$\int_U (k_2 - k_1) \leq \int_U |k_1| \leq \epsilon$.  We define
$h \in A(\Di)$ by the formula defining $k$ above, but with
$k_1$ replaced by $-k_2$.  Then $h$ is nonvanishing on
$\T$, and $\vert k \vert \leq \vert \phi \vert
\leq |h|^{-1}$ on $\T$.  Setting $r = k_2 - k_1$, we may write, for
fixed $z \in \Di$:
$$\int_{\T} \frac{e^{i \theta} + z}{e^{i \theta} - z} \; r(
e^{i \theta})
\; = \epsilon \int_{\T}
\frac{e^{i \theta} + z}{e^{i \theta} - z} \;  -  \;
 \epsilon \int_U \frac{e^{i \theta} + z}{e^{i \theta} - z}
\; +
\; \int_U \frac{e^{i \theta} + z}{e^{i \theta} - z}  \; r(e^{i \theta})
\;  \; . $$
The first of the three terms
 on the right equals $2 \pi \epsilon$.
Supposing that $d(z,U) \geq \delta$, we have
$$ |\frac{e^{i \theta} + z}{e^{i \theta} - z}| \leq \frac{2
}{\delta} \; \; ,
$$
for $e^{i \theta} \in U$, whence
$$| \frac{1}{2\pi} \int_{-\pi}^{\pi}
\frac{e^{i \theta} + z}{e^{i \theta} - z} \; (k_2-k_1)(e^{i
 \theta})
\; d \theta | \; \leq \; \epsilon (1 + \frac{4}{\delta}) \;
 .
$$

We now  check that (ii) of Definition
\ref{arr} holds.  By an easy compactness argument, we may assume
that the compact subset $C \subset \Omega \setminus E$ there,
is a finite closed interval.
Pick $\epsilon$ so small
in relation to $d(C,E)$, that for any $z$ close enough to
$C$, we have $d(z,s) \geq \sqrt{\epsilon}$ for any $s \in 
\T$
with $d(s,E) < \epsilon$.
For the $k_1, k_2$ associated
with this $\epsilon$, we have, for $z$ close to $C$, that
   $$ \vert \frac{1}{2\pi} \int_{-\pi}^{\pi}
\frac{e^{i \theta} + z}{e^{i \theta} - z} (k_2 - k_1)(e^{i \theta})
\; d \theta \vert \leq 5 \sqrt{\epsilon} \; \; .
$$
Thus $$\left| \exp \left( \frac{1}{2\pi} \int_{-\pi}^{\pi}
\frac{e^{i \theta} + z}{e^{i \theta} - z} (k_1 - k_2)(e^{i \theta})
\; d \theta \right)
\; \; - \; \; 1 \right| \leq 6 \sqrt{\epsilon} \; \; ,
$$
for all $z$ close enough to $C$ .
Thus $ |k(e^{i \theta})h(e^{i \theta}) - 1| \leq 6
 \sqrt{\epsilon} $,
on $C$.
From this it is clear that the conditions of  Definition
\ref{arr} are met, so that $A f$ is a rigged module over
$A(\Di)$.
\end{proof}

One may always choose the $\phi$ in the last theorem to be
an outer function, and then $|\phi|$ will also be
the `unique harmonic
extension' of $f$ to $M_A = \bar{\Di}$.

In this
case of $A = A(\Di)$, we see that any  $f \in \R$ is nonvanishing
inside $\Di$, and has a harmonic
logarithm on all of $\Di$.

\begin{proposition}
\label{sg}
The set $\R(\Omega)$ is a unital semigroup.
Indeed, if $f_1 \in \R(\Omega)$ and $f_2 \in C(\Omega)_+
$ then $(Af_1)^{\bar{}} \otimes_{hA} (Af_2)^{\bar{}}
\cong (Af_1 f_2)^{\bar{}}$ (completely) $A$-isometrically.
If $f_1 \in \R_{E_1}$, and if $f_2 \in \R_{E_2}$, then
$f_1 f_2 \in \R_{E_1 \cup E_2}$.
\end{proposition}

\begin{proof}  Clearly the multiplication map $\Phi :
(Af_1)^{\bar{}} \otimes_{hA} (Af_2)^{\bar{}}
\rightarrow (Af_1 f_2)^{\bar{}}$ is completely contractive,
and has dense range.
Conversely, choose $H_m, K_m$ as in Definition \ref{arr},
and let $e_m = H_m \cdot K_m$ as before.  For $a \in A$,
define $\theta_m(af_1 f_2)
 = e_m f_1 \otimes f_2 a = H_m f_1 \odot K_m f_2 a $.
The $\odot$ notation here, is commonly used with reference to
the Haagerup tensor product.  Namely, for two finite
tuples $x =(x_i), y = (y_i)$ , the expression
$x \odot y$ means $\sum_i x_i \otimes y_i$.
If $g = \chi_E f_1^{-1}$, then
$\theta_m(af_1f_2) = H_m f_1 \odot K_m g (af_1f_2)$.
From the definition of the Haagerup tensor product \cite{BMP}
we see from this
that $\theta_m$ is well defined and completely contractive.
It is easy to see that
 $\theta_m (\Phi(u)) \rightarrow u$ , for all
$u \in Af_1 \otimes_A Af_2$.  Hence $\Phi$ is a complete isometry.

We leave it to the reader to check that if
$f_1, f_2 \in \R(\Omega)$, then $f_1f_2 \in \R(\Omega)$.
Thus $\R(\Omega)$ is a unital semigroup.  The last assertion
is also an easy exercise.
\end{proof}

\vspace{5 mm}

Thus $\log \R_E$ and $\log
 \R(\Omega)$ are `harmonic classes' in the sense
of the introduction.

\section{Applications to Morita bimodules.}

\begin{theorem}
\label{whatname}  If $X$ is an algebraically singly generated faithful
function module (or equivalently, of the form $Ag$ for some $g \in 
C(\Omega)^+$), the following are equivalent:
\begin{itemize}
\item [(i)]  $X \cong Af$ $A$-isometrically, for some $f \in \M$;
\item [(ii)]  For any $\epsilon > 0$, there exists $n \in \N$,
and  $A$-module maps $\varphi : X \rightarrow A^{(n)}$ and 
$\psi : A^{(n)} \rightarrow X$, with $\psi \circ \varphi =
Id_X$, and $\Vert \varphi \Vert_{cb} \leq 1 + \epsilon$
and $\Vert \psi \Vert_{cb} = \Vert \psi \Vert \leq 1 + \epsilon$;
\item [(iii)]  $X$ is $A$-rigged.
\end{itemize}
\end{theorem}

In (ii), the operator space structure on $A^{(n)}$ is $R_n(A)$,
as in the definition of $A$-Hilbertian.

\begin{proof}  $(i) \Rightarrow (ii)$:  This follows from the
definition of $\M$; given $\epsilon > 0$ there exist $A$-tuples
$H, K$ with $H.K = 1$, such that the norms of 
$Hf$ and $Kf^{-1}$ are close to 1.  These may be associated 
with maps $\varphi$
and $\psi$ as in the proof of Corollary \ref{Hisr}.

$(ii) \Rightarrow (iii)$:  Follows from the definition of
$A$-Hilbertian, and the fact that such $X$ is $A$-rigged
if and only if it is $A$-Hilbertian.

$(iii) \Rightarrow (ii)$:   The maps $\phi_m, \psi_m$
in the definition of $A$-Hilbertian, may be associated as in
the proof of Corollary \ref{Hisr}, with certain $A$-tuples
$H_m, K_m$.  We have $(H_m.K_m) f = \psi_m(\phi_m(f)) \rightarrow
f$.  Hence $(H_m.K_m) \rightarrow 1$, so that by the common 
Neumann series trick, $H_m.K_m$ is an invertible element $b^{-1}$
of $A$, and $|b| \approx 1$.  Replace $K_m$ by $K_m b$.
Correspondingly we get an adjusted $\phi'_m, \psi'_m$, which now 
satisfies (ii).

$(ii) \implies (i)$ we state as the next result.
\end{proof}

\vspace{3 mm}

The following result is highly analagous to Theorem \ref{new}:

\begin{proposition}
\label{new3}   Suppose that
 $A$ is a uniform algebra on compact $\Omega$,
and that $f \in C(\Omega)^+$.  Then $X = Af$ satisfies condition 
(ii) (or equivalently (iii))
of the above theorem, if and only if $Af \cong (Af)_{|\partial A}$
isometrically via the restriction map, and $f_{|\partial A}
\in \M(\partial A)$.
\end{proposition}

\begin{proof}  The $(\Leftarrow)$ direction is ($(i) \Rightarrow (ii)$)
of the previous theorem.

Supposing (ii) of the Theorem, we proceed as in the proof of 
Corollary \ref{Hisr}, to associate with $\varphi$ and $\psi$,
$A$-tuples $H$ and $K$.  We have that $1=H.K$, and as in that 
proof we get $f(\cdot) \Vert H(\cdot) \Vert_2 \leq 1+\epsilon$
on $\Omega$, and $\Vert K(\cdot) \Vert_2 \leq
(1+\epsilon) f(\cdot)$ on $\partial A$.   Applying
Cauchy-Schwarz to $1 = |(Hf).(Kf^{-1})|$, and using 
these inequalities, shows that 
$(1+\epsilon) \Vert K(\cdot) \Vert_2 \geq f$ on $\Omega$,
and 
$\Vert H(\cdot) \Vert_2 \approx f^{-1}$ on $\partial A$. 
Thus $f_{|\partial A}
\in \M(\partial A)$.   That $Af \cong (Af)_{|\partial A}$
follows as in the proof of \ref{Gmax}.
\end{proof}

\begin{corollary}
\label{wenmo}  Suppose that 
$X$ is a {\em singly generated} left
$A$-rigged module.  Then the following are equivalent:
\begin{itemize}
\item [(i)]   The peak set $E$ associated with $X$ is the empty set,
\item [(ii)]  $X$, with the obvious (symmetric) right module
action,  is a strong Morita equivalence
$A-A$-bimodule.
\item [(iii)]  $X$ is algebraically singly generated.
 \end{itemize} 
If $M_A$ is connected, 
then the above are also equivalent
to:
\begin{itemize}
\item [(iv)]  $X$ is {\em algebraically} finitely generated
and projective as a left $A$-module.
\end{itemize}
\end{corollary}

\begin{proof}  $(i) \Rightarrow (ii)$:
Follows since $\M = \R_{\emptyset}$.

$(ii) \Rightarrow (i)$ :  Assuming (ii), then by Theorem \ref{nice},
$X \cong Af$ $A$-isometrically, where $f \in C(\partial A)^+$.
Since every strong Morita equivalence bimodule is a rigged module,
Corollary \ref{H} (iii) now gives (i).

$(ii) \Rightarrow (iv)$: This is  true  for 
any algebraic Morita equivalence
$A-A$-bimodule \cite{F}.
 
$(iv) \Rightarrow (i)$: It follows by 
\cite{Bhmo} Theorem 3.6 (6), and \ref{just} below,
that $\K(X) = J_E$ is unital.   Thus $\chi_E$ is continuous, so that 
$E$ is closed and open.  If $\Omega$ is connected, 
it follows that $E$ is the empty set.

$(iii) \Rightarrow (ii)$:  Every singly generated  $A$-rigged 
module is of the form $(Af)^{\bar{}}$, which is a faithful 
function $A$-module.   
Now appeal to ($(iii) \Rightarrow (i)$) of
Theorem \ref{whatname}.
\end{proof}

The proof of $(ii) \Rightarrow (i) \Rightarrow (ii)$
above, together with the basic theory of strong Morita equivalence
(see \cite{BMP} section 4), 
 gives a strengthened form of Theorem \ref{nice}:   

\begin{corollary}
\label{Motat}  Let $X$ be a singly generated 
symmetric
strong Morita equivalence $A-A$-bimodule.  Then there
exists an $f \in \M$ such that $X \cong Af$ $A$-isometrically,
$\tilde{X} \cong B_A(X,A) \cong Af^{-1}$ $A$-isometrically,
and the strong  Morita context associated with $X$ may be 
identified with $(A,A,Af,Af^{-1})$.
\end{corollary}

Finally, we show that every singly generated rigged module over
a function algebra, is a strong Morita equivalence bimodule
(over a possibly different algebra):
 
\begin{theorem}
\label{just}
\begin{itemize}
\item [(i)]
If $X$ is a singly generated operator module over a
function algebra $A$ on a compact space $\Omega$, then
$X$ is an $A$-rigged module if and only if
there exists a peak set $E$ in $\Omega$ for $A$, and a
function $f \in \R_E(\Omega)$,
such that  $X \cong (J_E f)^{\bar{}}$
$A$-isometrically.
\item [(ii)]  If the equivalent conditions in (i) hold, then
 $X$ , with the obvious (symmetric) right module action,
is a strong Morita equivalence $J_E-J_E$-bimodule.
\item [(iii)] Conversely, if $E$ is a p-set, then
 any strong Morita equivalence
$J_E-J_E$-bimodule,
or more generally any $J_E$-rigged module,
is an $A$-rigged module.
\end{itemize}
\end{theorem}

\begin{proof}   The proof requires some technical knowledge of rigged modules
\cite{Bhmo}.

 Here is
one way to see (iii).
If $X$ is a
left $J_E$-rigged module, then from
Corollary \ref{co} and \S 6 of \cite{Bhmo},
 $J_E \otimes_{hJ_E} X$ is
a left $A$-rigged module.  But $J_E \otimes_{hJ_E} X \cong X$, $A$-isometrically.

To get (i) and (ii), suppose that $X$ is a singly generated
$A$-rigged module.  Then by the previous results, there is
a peak set $E$ and a function $f \in \R_E(\Omega)$ such that
$X \cong (Af)^{\bar{}}$ $A$-isometrically.
Since $f$ is a strictly positive element of
the ideal $I$ in Lemma \ref{rig}, we see that
$I = \{ p \in C(\Omega) :
p(x) = 0 \; \text{for all} \; x \in E \}$.
Suppose that $Y = \tilde{X}$ is the dual rigged module of
$X$ (see \cite{Bhmo,BMP} for details),
where the reader may also
find the definition of $\K(X)$, which we shall need shortly.
Since $W$ may be taken to be
$I$, it follows that
the linking C$^*-$algebra
for $W$ is $M_2(I)$.
Thus we can make
the following deductions from Theorem 5.10 in \cite{Bhmo}.
Firstly, $\K(X)$ may be identified with a
closed subalgebra $J$ of $I$,
and $J$ has a contractive approximate identity which is a
contractive approximate identity for $I$.  Also, $Y$ may be
regarded as a subspace of $I$, and the canonical pairings
$X \times Y \rightarrow A $, and $Y \times X \rightarrow
\K(X) = J$ may be regarded as the commutative
multiplication in $I$ .  Thus it follows that
$J$ is the
closure of the span of the range of the
canonical pairing $X \times Y \rightarrow
A$ .  Thus $J$ is also a subset, indeed a
closed ideal, of
$A$.   By the commutativity of the multiplication in $I$,
$X$ is a strong Morita equivalence
$J-J$-bimodule.  We also see that $JX$ is dense in $X$.

Since $J$ has a
contractive approximate identity, we may appeal to Theorem
\ref{Smith} to see that $J = J_{E'}$ for a p-set $E'$.
Since $J$ contains a c.a.i. for $I$, it is
clear that $E' \subset E$.  On the other hand, if
$w \in E \setminus E'$, then there is a peak set $F$ containing
$E'$, with $w \notin F$.  Choose $a \in A$ with $a \equiv 0$
on $E'$ and $a(w) \neq 0$.  Then $a \in J$, so $a \in I$.
This is impossible, so that $E = E'$.

Note that $Af$ is dense in $X$, $JX$ is dense in $X$ and $JA = J$.
Hence $Jf$ is dense in $X$.  Thus $X = (Jf)^{\bar{}}$.
 \end{proof}

\vspace{5 mm}

We thank C. Le Merdy, R. R. Smith, and V. Paulsen
for helpful conversations.

\end{document}